\documentclass[
,twocolumn%
,tightenlines%
,amssymb,amsmath,nobibnotes,aip,chaos,showpacs]{revtex4}
\usepackage{amsmath,amsfonts,amssymb,amscd,amsthm,epsf}
\usepackage{rotating}
\usepackage{comment}

\usepackage{bm}% bold math
\usepackage{amsthm}
\usepackage{amssymb}
\usepackage{amsfonts}
\usepackage{mathrsfs}
\usepackage{amsbsy}

\usepackage{dcolumn}% Align table columns on decimal point

\usepackage{comment}

\usepackage{epsf}
\usepackage{epsfig}
\usepackage{graphicx}

\newcommand{\calD}{{\mathcal{D}}}

\newcommand{\calF}{{\mathcal{F}}}

\newcommand{\bbC}{{\mathbb{C}}}

\newcommand{\bbN}{{\mathbb{N}}}

\newcommand{\bbT}{{\mathbb{T}}}

\newcommand{\sfA}{{\mathsf{A}}}
\newcommand{\sfB}{{\mathsf{B}}}

\newcommand{\rme}{{\mathrm{e}}}

\newcommand{\rmS}{{\mathrm{S}}}

\def\complex{{\mathbb C}}

\newcommand{\siegeldisk}{\calD}
\newcommand{\ds}[1]{\displaystyle{#1}}

\newcommand{\regul}{{\kappa}}

\newtheorem{theorem}{Theorem}[section]

\newtheorem{preremark}[theorem]{Remark}

\begin{document}

%\preprint{}

\title{Boundaries of  Siegel disks -- numerical studies 
        of their dynamics and regularity}


\author{Rafael de la Llave}
\email[]{llave@math.utexas.edu}
%\thanks{}
\affiliation{Department of Mathematics, University of Texas, 
        Austin, TX 78712, USA}

\author{Nikola P.\ Petrov}
\email[]{npetrov@math.ou.edu}
\affiliation{Department of Mathematics, University of Oklahoma, 
        Norman, OK 73019, USA}


\date{May 12, 2008}

\begin{abstract}
Siegel disks are domains around fixed points of holomorphic maps 
in which the maps are locally linearizable 
(i.e., become a  rotation under an appropriate change of 
coordinates which is analytic in a neighborhood of the origin). 
The dynamical behavior of the iterates of the map 
on the boundary of the Siegel disk 
exhibits strong scaling properties 
which have been intensively studied 
in the physical and mathematical literature.  

In the cases we study, 
the boundary of the Siegel disk is a Jordan curve 
containing a critical point of the map 
(we consider critical maps of different orders), 
and there exists a natural parameterization 
which transforms the dynamics on the boundary into a rotation.  
We compute numerically this parameterization and 
use methods of harmonic analysis 
to compute the global H\"older regularity 
of the parameterization 
for different maps and rotation numbers.

We obtain that the regularity of the boundaries 
and the scaling exponents 
are universal numbers in the sense of renormalization theory
(i.e., they do not depend on the map when the map
ranges in an open set), 
and only depend on the order of the critical point 
of the map in the boundary of the Siegel disk 
and the tail of the continued function expansion 
of the rotation number. 
We also discuss some possible relations 
between the regularity of the parameterization 
of the boundaries and the corresponding 
scaling exponents.  
\end{abstract}


%      05.45.Df, %      Fractals (see also 47.53.+n Fractals in fluid dynamics
%      02.30.Nw, %      Fourier analysis
%      05.10.-a, %      Computational methods in statistical physics 
%               %      and nonlinear dynamics (see also 02.70.-c in 
%               %      mathematical methods in physics)
%      05.10.Cc %      Renormalization group methods
       %05.45.Pq %numerical simulations,

%\keywords{????????????}

\maketitle

%\section{}
%\subsection{}
%\subsubsection{}

%\tableofcontents

{\bf 
According to a celebrated theorem by Siegel, 
under certain arithmetic conditions, 
the dynamical behavior of the iterates 
of a holomorphic map around 
a fixed point of the map is very simple 
-- the iterates of the map fill densely 
analytic  topological circles around the critical point.  
In the domain around the critical behavior 
where the iterates exhibit such behavior 
-- called Siegel disks, 
-- there exists a complex analytic  change of variables 
that makes the map locally a multiplication 
by a complex number of modulus~1.  
On the boundary of a Siegel disk, however, 
the dynamical behavior of the iterates of the map 
is dramatically different 
-- for example, the boundary is not a smooth curve.
The iterates on the boundary of the Siegel disk 
exhibit scaling properties 
that have motivated the development 
of a renormalization group description.  
The dynamically natural parameterization 
of the boundary of the Siegel disk has low regularity.  
We compute accurately the natural parameterization of 
the boundary 
and apply methods from Harmonic Analysis 
to compute the H\"older exponents 
of the parameterizations of the boundaries 
for different maps 
with different rotation numbers and with different 
orders of criticality of their critical points.  
%Although the boundaries of Siegel disks 
%do not have direct physical meaning, 
%they are a model for a critical transition.  
}

\section{Introduction}

Siegel disks -- the domains around fixed points 
of holomorphic maps in which the map is locally 
linearizable (defined in more detail 
in Section~\ref{sec:disks-bound}) 
-- are among the main objects of interest in 
the dynamics of holomorphic maps. 
Their boundaries have surprising geometric properties 
which have attracted the attention of both mathematicians 
and physicists.  Notably, 
it was discovered in \cite{MantonN83, Widom83} 
that in some cases there were scaling relations for 
the orbit, which suggested that the boundary was a fractal object.  
Since then this phenomenon has been a subject 
of extensive numerical and mathematical studies 
\cite{Osbaldestin92,Stirnemann94b,BurbanksS95,BishopJ97,McMullen98,BurbanksOS98-EPJB,BurbanksOS98,GraczykJ02,GaidashevY06,Gaidashev07}.  

In this paper, we report some direct numerical calculations of 
the H\"older regularity of these boundaries 
for different rotation numbers of bounded type
and for different maps. 

The main conclusion of the numerical calculations 
in this paper, is that, for the cases we consider, 
the boundaries of the Siegel 
disks are $C^\regul$ curves for some $\regul > 0$, 
and we can compute numerically the value of $\regul$. 
Even if we -- obviously -- consider only a finite number of 
cases, we expect that the results are significative for 
the Siegel disks of polynomials with rotation numbers which 
have an eventually periodic continued fraction. 

The values of the H\"older regularity $\regul$ are, up to the error of 
our computation, {\em universal} in the sense of 
renormalization group analysis, namely that they are independent 
of the map in a small neighborhood in the space of maps.  
We also performed computations for maps whose rotation numbers 
have the same ``tail'' of their continued fraction 
expansion, and found that our numerical results 
depend only on the tail.  

Our computation of the H\"older regularity 
are based on the method introduced in \cite{LlaveP02}, 
which is a numerical implementation of 
several constructions in Littlewood-Paley theory. 
This method was also used in \cite{Carletti03,ApteLP05,FuchssWAM06,OlveraP06}.

For the case of the golden mean rotation number, 
the fact that the boundaries of 
Siegel disks are H\"older  was proved in \cite{BurbanksS95} assuming the 
existence of the fixed point of the renormalization operator 
conjectured in~\cite{Widom83}. 
The existence of a  fixed point of a slightly different 
(and presumably equivalent) 
renormalization operator 
was proved in~\cite{Stirnemann94b}. It seems that a similar 
argument  will work for other rotation numbers 
with periodic continued fraction expansion
 provided that one has a fixed point of the appropriate 
renormalization operator.
 These arguments provide  
bounds to the H\"older regularity~$\regul$, based 
on properties of the fixed point of the renormalization operator.

Our computations rely on several rigorous mathematical results 
in complex dynamics.  
Notably, we will use that 
that for bounded type rotatation numbers 
and polynomial maps, 
the boundary of the Siegel disk is a Jordan curve, 
and contains a critical point \cite{GraczykS03}.
See also \cite{Ghys84,Rogers98} for other results in 
this direction.
We recall that a number is bounded type means 
that the entries in the continued fraction of this number 
are bounded. Equivalently, a number $\sigma$ is of
constant type if and only if for every natural
$n$ and integer $m$ we have $|\sigma n - m|^{-1} \le \nu n^{-1}$
\cite{Herman79}. 

It is also known in the mathematical literature 
that for non-Diophantine rotation numbers 
(a case we do not consider here and which indeed seems 
out of reach of present numerical experiments), 
it is possible to make the circle smooth 
\cite{Perezmarco99,AvilaBC04,BuffC07} 
or, on the contrary, not a Jordan curve~\cite{Herman85c, Rogers95}.  

The plan of this paper is the following. 
In Section~\ref{sec:disks-bound} 
we give some background on Siegel disks and explain 
how to parameterize their boundaries, 
Section~\ref{sec:num-regul} is devoted 
to the numerical methods used and the results 
on the regularity of the boundaries of the Siegel disks 
and the scaling properties of the iterates.  
In Section \ref{sec:radius-area} we consider 
some connections with geometric characteristics 
of the Siegel disk (in particular, its area), 
in Section \ref{sec:bound-regularity} we derive 
an upper bound on the regularity, 
and in the final Section~\ref{sec:discussion} 
we recapitulate our results.  
%and give possible directions of further study.  

%%%%%%%%%%%%%%%%%%%%%%%%%%%%%%%%%%%%%%%%%%%%%%%%%%%%%%%%%%%%
%%%%%%%%%%%%%%%%%%%%%%%%%%%%%%%%%%%%%%%%%%%%%%%%%%%%%%%%%%%%
%%%%%%%%%%%%%%%%%%%%%%%%%%%%%%%%%%%%%%%%%%%%%%%%%%%%%%%%%%%%

\section{Siegel disks and their boundaries\label{sec:disks-bound}}

\subsection{Some results from complex dynamics\label{sec:prelim-complex}}

In this section we summarize some facts from complex dynamics, 
referring the reader to \cite{Milnor06,Lyubich86} for more details.  

We consider holomorphic maps of $\complex$
that have a fixed point, 
and study their behavior around this point. 
Without loss of generality, in this section 
we assume that the fixed point is 
the origin, so that the maps have the form 
\begin{equation}  \label{eq:f-general}
f(z)  = a z + O(z^2) \ .
\end{equation}  
(For numerical purposes, we may find more efficient 
to use another normalization.)  
We are interested in the case that $|a| = 1$, 
i.e., $a = \rme^{2 \pi i \sigma}$, 
where $\sigma \in [0,1)$ is called the {\em rotation number} of~$f$.  
In our case, we take $f$ to be a polynomial, so that 
the domain of definition of the map $f$ is not an issue.  

The stability properties of the fixed point depend crucially 
on the arithmetic properties of~$\sigma$. 
The celebrated {\em Siegel's Theorem} 
\cite{Siegel42,Moser66a,Zehnder77} 
guarantees that, if $\sigma$ satisfies 
some arithmetic properties ({\em Diophantine conditions}), 
then there is a unique  analytic mapping $h$ 
(called ``conjugacy'') 
from an open disk of radius $r$ (called the 
Siegel radius)  around the origin, $B(0,r)$, 
to $\complex$ in such a way 
that $h(0) = 0$, $h'(0) = 1$, and 
\begin{equation}\label{conjugacy} 
f\circ h(z) = h(a z) \ .
\end{equation} 
We note that the Siegel radius is a 
geometric property of the Siegel disk. 
It is shown in \cite{SiegelM95} that $h$ 
can be characterized as the conformal mapping from 
$B(0,r_s)$ to the Siegel disk mapping $0$ to $0$ 
and having derivative $1$. Later in 
Section~\ref{sec:radius}, we will show how 
the Siegel disk can be computed effectively in 
the cases we consider.

We refer to \cite{Herman86,Douady87} for some mathematical developments 
on improving the arithmetic conditions of the Siegel theorem.
In this paper we only consider rotation numbers that satisfy 
the strongest possible Diophantine properties.  
Namely, we assume that $\sigma$ is of {\em bounded type} 
and, in particular, perform our computations 
for numbers with eventually periodic continued fraction expansions 
(see Section~\ref{sec:cont-frac} for definitions).  
In this case, there is an elementary proof 
of Siegel's Theorem~\cite{Llave83}.  

Let $r_\rmS$ stand for the radius of the largest disk 
for which the map $h$ exists.  
The image under $h$ of the open disk $B(0,r_\rmS)$ 
is called the {\em Siegel disk}, $\siegeldisk$, of the map~$f$.  
For $r<r_\rmS$, the image of each circle $\{w\in\complex:|w|=r\}$ 
under $h$ is an analytic circle.  
The boundary, $\partial \siegeldisk$, of the Siegel disk, 
however, is not a smooth curve for the cases 
considered here.  

%Boundaries of Siegel disks have attracted attention from 
%mathematicians and physicists. 
The paper \cite{MantonN83} 
contains numerical observations that suggest 
that the dynamics of the map $f$ on $\partial \siegeldisk$ 
satisfies some scaling properties. 
These scaling properties were explained in certain cases 
by renormalization group analysis 
\cite{Widom83,Stirnemann94a,Stirnemann94b, GaidashevY06,Gaidashev07}. 
These scaling properties suggest that the boundaries of 
Siegel disks can be very interesting fractal objects. 

Clearly, the Siegel disk cannot contain critical points 
of the map~$f$.  
It was conjectured in \cite{MantonN83} 
that the boundary of the Siegel disk contains a critical point.  
The existence of critical points on the boundary depends on the arithmetic 
properties of $\sigma$ and it may be false \cite{Herman85c}, 
but it is true under the condition that the rotation number is 
of bounded type \cite{GraczykS03}, which is the case 
we consider in this paper. See also~\cite{Rogers98a}. 

\begin{figure}
\centering
        \includegraphics[width=0.5\textwidth]{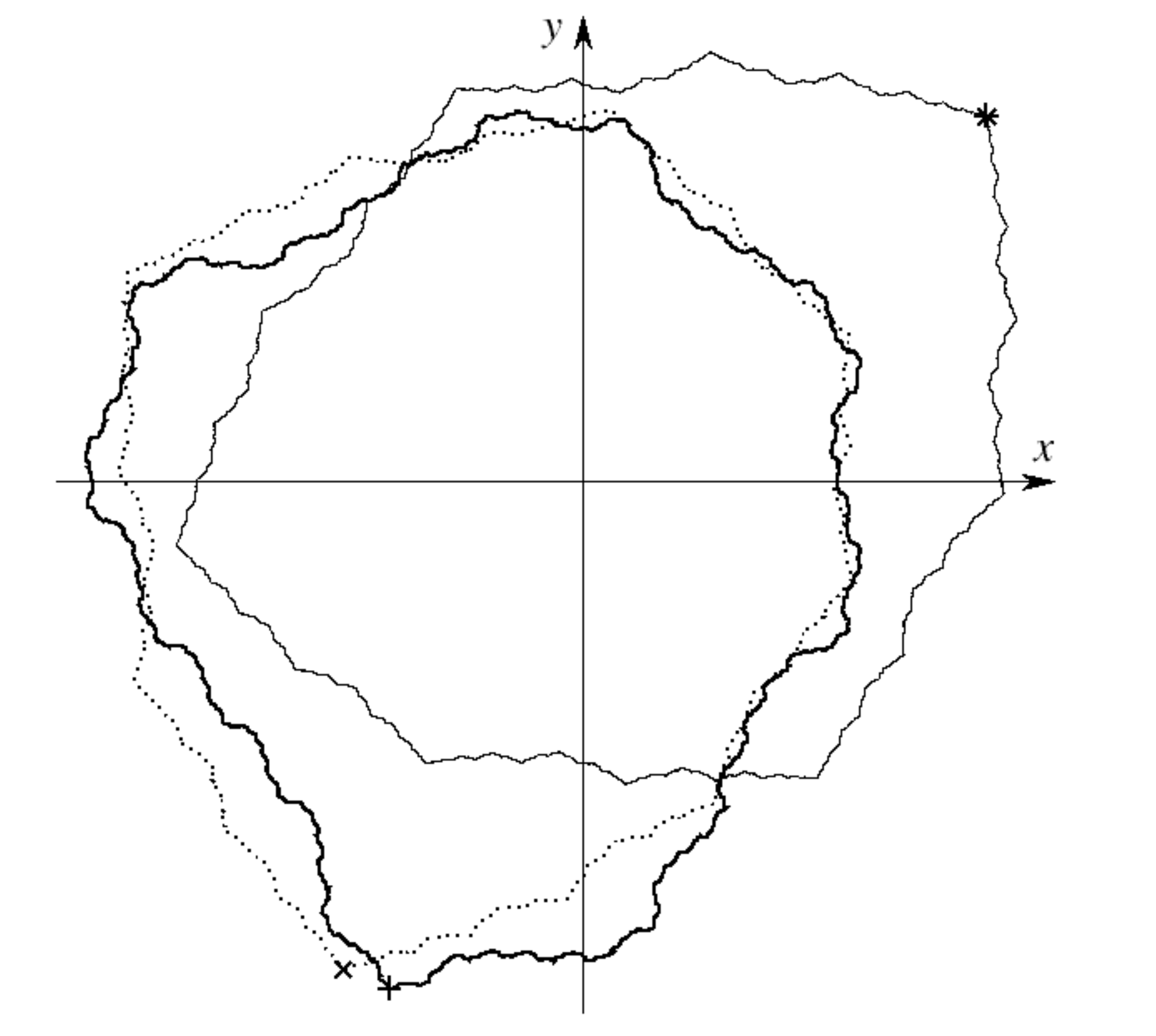}%width=0.55\textwidth
        \caption{\label{fig:orbits_1_inf_and_5_inf_and_5_1_inf}
        Critical points $c$ and Siegel disks of the map \eqref{eq:old-f}
        for $\sigma=\langle 1^\infty\rangle=\frac12(\sqrt{5}-1)$ 
        ($*$ and a thin line), 
        $\sigma=\langle 5^\infty\rangle=\frac12(\sqrt{29}-5)$ 
        ($+$ and a thick line), 
        and 
        $\sigma=\langle 51^\infty\rangle=\frac1{38}(9-\sqrt{5})$ 
        ($\times$ and a dotted line).  
        }
\end{figure}

In Figure~\ref{fig:orbits_1_inf_and_5_inf_and_5_1_inf} 
we show the Siegel disks of the map 
\begin{equation}  \label{eq:old-f}
f(z)=\rme^{2\pi i\sigma}z+z^2 
\end{equation}
for different rotation numbers~$\sigma$ 
(for the notations for $\sigma$ 
see Section~\ref{sec:cont-frac}).  
In all cases the only critical point, 
$c=-\frac12\rme^{2\pi i\sigma}$, 
is simple: 
$f'(c)=0$, $f''(c)=2\neq 0$.  

In this paper we study maps of the form 
\begin{equation}  \label{eq:f-def}
f_{m,\sigma,\beta} (z) 
= \frac{1}{\beta} \rme^{2\pi i\sigma} 
        \left[g_{m+1}(z) - (1-\beta) g_{m}(z) \right] \ ,
\end{equation}
where $m\in\bbN$, $\beta$ is a complex parameter, 
and the function $g_m:\bbC\to\bbC$ is defined as 
\[
g_m(z) = \frac{1}{m+1} \left[1-(1-z)^{m+1}\right] \ .
\]

Let $f$ be a map of the form \eqref{eq:f-general}, 
and $c$ be its critical point 
that belongs to the boundary of the Siegel disk 
of this map (we will only consider cases 
where $\partial\calD$ contains one critical point).  
Let $d$ be the {\em multiplicity} 
of the critical point $c$, 
i.e., $f^{(k)}(c)=0$ for $k=1,2,\ldots,d$, 
and $f^{(d+1)}(c)\neq 0$.  
We will call $d$ the {\em order} of the critical point.  

Noticing that, for the map \eqref{eq:f-def}, 
\[
f'_{m,\sigma,\beta}(z) 
= 
\rme^{2\pi i\sigma} (1-z)^m \left(1-\frac{z}{\beta}\right) \ ,
\]
we see that $f_{m,\sigma,\beta}(0)=0$, 
$f'_{m,\sigma,\beta}(0)=\rme^{2\pi i\sigma}$, 
and, more importantly, if $\beta\neq 1$, 
the point $z=1$ is a zero of 
$f'_{m,\sigma,\beta}$ of multiplicity~$m$, 
while, for $\beta=1$, 
the point $z=1$ is a zero of 
$f'_{m,\sigma,\beta}$ of multiplicity~$m+1$.  
%We will call the multiplicity $d$ 
%the {\em order} of the critical point of the map.  
As long as the critical point $z=\beta$ 
is outside the closure of the Siegel disk, 
the scaling properties of the iterates 
on $\partial\calD$ in the case of Diophantine $\sigma$ 
are determined by the order $d$ of the critical point $z=1$.  
Below by ``critical point'' we will mean 
the critical point that belongs to $\partial\calD$.  
One of the goals of this paper is to study 
how the regularity and the scaling properties 
depend on the order of the critical point.  

\begin{figure}
\centering
        \includegraphics[width=0.49\textwidth]{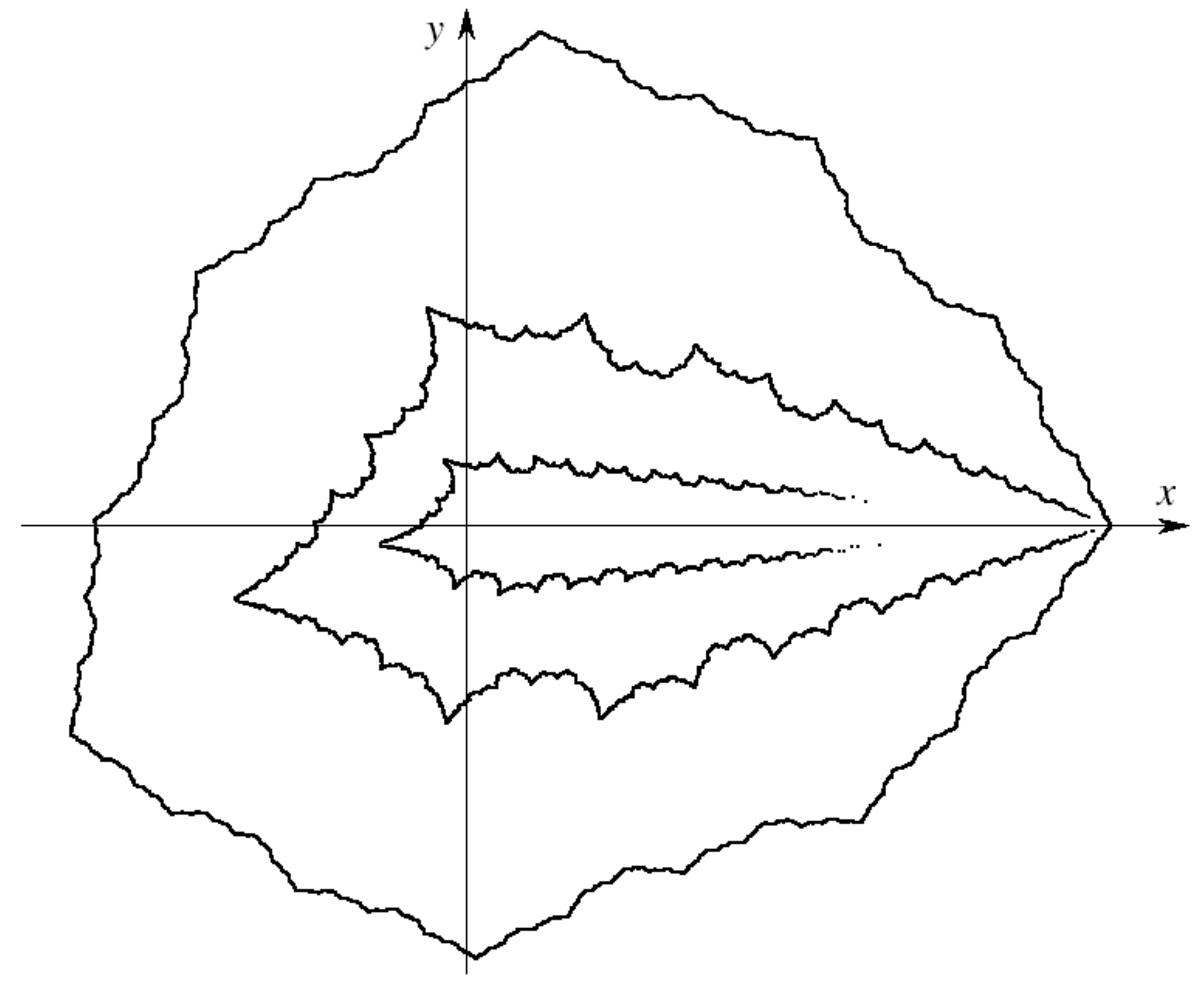}%width=0.49\textwidth
        \caption{\label{fig:high-crit}
        Siegel disks of the maps 
        $f_{d,\langle 2^\infty\rangle,1+3i}$ 
        (given by \eqref{eq:f-def}) 
        with critical point $c=1$ of order 
        $d=1$ (outermost curve), 
        $d=5$, 
        and $d=20$ (innermost curve) 
        -- see text for details.  
        }
\end{figure}

In Figure~\ref{fig:high-crit} we show about 16 million iterates 
of the critical point $c=1$ of the maps 
$f_{d,\langle 2^\infty\rangle,1+3i}$ 
for order $d=1$, 5, and 20, of the critical point $c=1$ 
(the other critical point, $\beta=1+3i$, 
is not in the closure of the Siegel disks, 
so is irrelevant for the problem studied).  
Note that, especially for highly critical maps, 
the iterates approach the critical point very slowly 
because the modulus of the scaling exponent 
becomes close to~1 
(see Table~\ref{table:exponents}).

%%%%%%%%%%%%%%%%%%%%%%%%%%%%%%%%%%%%%%%%%%%%%%%%%%%%%%%%%%%%

\subsection{Parameterization of the boundary 
        of a Siegel disk\label{sec:parameterization}}

In the cases considered here, 
the boundaries of Siegel disks cannot be written 
in polar coordinates as $r=R(\theta)$ 
(because some rays $\theta=\mathrm{const}$ 
intersect $\partial\siegeldisk$ more than once).  
In this section 
we explain how to parameterize $\partial\siegeldisk$, 
and define the functions whose regularity 
we study numerically.  

Once we know that a critical point $c$ is in the boundary of 
the Siegel disk $\siegeldisk$ (which in the cases 
we consider is guaranteed by the 
results of \cite{GraczykS03}), 
it is easy to obtain a parameterization of the 
boundary which semiconjugates the map $f$ to a rotation.  

It is known from the mathematical theory that $h$ 
-- which is univalent in the open disk $B(0,r_\rmS)$ 
-- can be extended to the boundary of $B(0,r_\rmS)$ 
as a continuous function 
thanks to the Osgood-Taylor-Carath\'eodory Theorem 
(see, e.g., 
\cite[Section 16.3]{Henrici-ACCA-1} or 
\cite[Section IX.4]{Burckel-I}).  

A dynamically natural parameterization 
$\chi$ of $\partial\siegeldisk$ 
is obtained by setting 
\begin{equation}\label{particular}
\chi(t) = h\left(r_\rmS \, \rme^{2\pi i(t+\theta)}\right) \ ,
\end{equation}
where $\theta$ is a constant to be chosen later.  
From \eqref{conjugacy}, it follows that 
\begin{equation}\label{conjugacy2}
f\circ \chi(t)  =  \chi(t + \sigma) 
\end{equation}
with $\chi(0) = c$.  
Since we know that in our cases the critical point $c$ 
is in $\partial\siegeldisk$, 
$|h^{-1}(c)| = r_\rmS$, 
and we can choose $\theta$ so that 
\begin{equation}  \label{chi0}
\chi(0) = h(r_\rmS\,\rme^{2\pi i \theta}) = c \ .
\end{equation}
In summary, 
the function $\chi:\bbT\to\bbC$ defined 
by \eqref{particular} and \eqref{chi0} is a parameterization 
of $\partial\siegeldisk$ such that 
the dynamics of the map $f$ on $\partial\siegeldisk$ 
is a rotation by $\sigma$ on the circle $\bbT$, 
and $\chi(0)$ is the critical point. 

Iterating $n$ times \eqref{conjugacy2} 
for $t=0$, and using \eqref{chi0}, we obtain 
\begin{equation}  \label{iters}
\chi(n\sigma) = f^n(c) \ .
\end{equation}
Since the rotation number $\sigma$ of $f$ is irrational, 
the numbers $n\sigma$ (taken mod~1) 
are dense on the circle $\bbT$, 
and the iterates $f^n(c)$ of the critical point $c$ 
are dense on $\partial\siegeldisk$ as well.  
Hence, using equation \eqref{iters}, 
we can compute a large number of values 
of $\chi$ by simply iterating~$f$.  

\begin{figure}
        \includegraphics[width=0.5\textwidth]{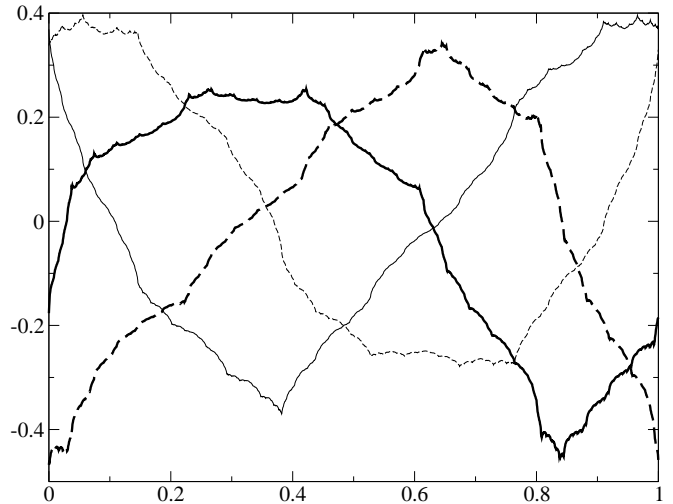}%width=0.49\textwidth
        \caption{\label{fig:functions1}
        Graphs of the real (solid lines) 
        and the imaginary (dashed lines) 
        parts of the map $\chi$ corresponding 
        to $\partial\siegeldisk$ of the map \eqref{eq:old-f} 
        for $\sigma=\langle 1^\infty\rangle$ 
        (thin lines) 
        and $\sigma=\langle 5^\infty\rangle$ 
        (thick lines) 
        -- cf.\ Figure~\ref{fig:orbits_1_inf_and_5_inf_and_5_1_inf}.  
        }
\end{figure}

Our main object of interest 
are the real, $\Re\chi$, and imaginary, $\Im\chi$, parts 
of the function~$\chi$.  
%In Figure~\ref{fig:real_imag_1_inf_5_inf} 
In Figure~\ref{fig:functions1} 
we have shown the graphs of the real and imaginary 
parts of $\chi$ parameterizing the boundary 
of the Siegel disk of the map \eqref{eq:old-f} 
for two different rotation numbers.  

%%%%%%%%%%%%%%%%%%%%%%%%%%%%%%%%%%%%%%%%%%%%%%%%%%%%%%%%%%%%
%%%%%%%%%%%%%%%%%%%%%%%%%%%%%%%%%%%%%%%%%%%%%%%%%%%%%%%%%%%%
%%%%%%%%%%%%%%%%%%%%%%%%%%%%%%%%%%%%%%%%%%%%%%%%%%%%%%%%%%%%

\subsection{Continued fraction expansions 
        and rational approximations\label{sec:cont-frac}}

In this section, we collect some of the notation on 
continued fraction expansions.

Let $\sfA=(a_1,a_2,\ldots,a_p)$ 
be a finite  sequence of $p$ natural numbers $a_j\in \bbN$; 
for brevity, we will usually omit the commas and write $\sfA=(a_1a_2\ldots a_p)$.  
Let $\sfB=(b_1b_2\ldots b_q)$ be another 
sequence of natural numbers $b_j\in\bbN$, $j=1,2,\ldots,q$, 
$\sfA\sfB:=(a_1 a_2 \ldots a_p b_1 b_2 \ldots b_q)$ 
stand for the concatenation of $\sfA$ and~$\sfB$, 
and $\sfA^n$ stand for $\sfA\sfA\cdots\sfA$ ($n$ times). 
Let $|\sfA|=p$ denote the length of~$\sfA$.  

For $a\in\bbN$ define the function 
$\calF_a:(0,1)\to(0,1)$ by $\calF_a(x) := \ds{\frac{1}{a+x}}$.  
Similarly, for $\sfA=(a_1a_2\ldots a_p)$, define the function 
$\calF_\sfA:(0,1)\to(0,1)$ as the composition 
\[
\calF_\sfA(x) := \calF_{a_1}\circ\calF_{a_2}\circ\cdots\circ\calF_{a_p}(x) \ .
\]
Let $\langle \sfB \rangle$ stand for the number whose continued 
fraction expansion (CFE) is given by the numbers in the sequence~$\sfB$:  
\[
\langle \sfB \rangle 
= 
\langle b_1 b_2 \ldots b_q \rangle 
:= 
\ds{\frac{1}{b_1+\ds{\frac{1}{b_2+\ds{\frac{1}{\ddots +\ds{\frac{1}{b_q}}}}}}}} \ ; 
\]
the numbers $b_j$ are called the {\em (partial) quotients} of~$\langle \sfB\rangle$.  
A number $\sigma=\langle a_1 a_2 \ldots \rangle$ 
is of {\em bounded type} if all numbers $a_i$ ($i\in\bbN$) 
are bounded above by some constant~$M$.  

We are especially interested in studying numbers with CFEs of the form 
\[
\langle \sfA\sfB^\infty \rangle
:= 
\lim_{n\to\infty} \langle \sfA\sfB^n \rangle \ ,
\]
which are called {\em eventually periodic} (or {\em preperiodic}).  
Since each number of this type is a root of a quadratic equation 
with integer coefficients (see \cite[Theorem 176]{HardyW-theory-numbers}), 
such numbers are also called {\em quadratic irrationals}.  
We will call $\sfA$ the {\em head} and $\sfB^\infty$ the {\em tail}, 
$\sfB$ the {\em period}, and $|\sfB|$ the {\em length} of the period of the CFE.  

If two quadratic irrationals have the same tail, they are said to be 
{\em equivalent}. 
 One can prove that $\sigma$ and $\rho$ are equivalent 
if and only if $\sigma=(\iota\rho+\lambda)/(\mu\rho+\nu)$, 
where the integers $\iota$, $\lambda$, $\mu$ and $\nu$ 
satisfy $\iota\nu-\lambda\mu=\pm 1$ 
\cite[Theorem 175]{HardyW-theory-numbers}.

%%%%%%%%%%%%%%%%%%%%%%%%%%%%%%%%%%%%%%%%%%%%%%%%%%%%%%%%%%%%

\subsection{Scaling exponents \label{sec:scaling-exp}}

Let $f$ be a map of the form \eqref{eq:f-general} 
with an eventually periodic 
rotation number $\sigma=\langle\sfA\sfB^\infty\rangle$ 
with length of its period $q=|\sfB|$, 
and let $c$ be the critical point of $f$ 
on~$\partial\siegeldisk$ 
(and there are no other critical points on~$\partial\siegeldisk$).  
Let 
\begin{equation}  \label{eq:fibonacci}
\frac{P_m}{Q_m}=\langle \sfA \sfB^m \rangle \ ,
\end{equation}
where $P_m$ and $Q_m$ are natural numbers 
that have no common factors.  
Define the {\em scaling exponent} 
\begin{equation}  \label{eq:scaling-exp}
\alpha 
:= 
\lim_{m\to\infty} 
\ds{\frac{f^{Q_{m+1}}(c)-c}{f^{Q_{m}}(c)-c}}  \ .
\end{equation}
This exponent is a complex number that depends 
on the tail $\sfB$ of the CFE of $\sigma$ 
and on the order of the critical point~$c$, 
but does not depend on the head~$\sfA$ 
or on details about the map~$f$.  

If the length $|\sfB|$ of the tail $\sfB$ of the CFE 
of the rotation number of the map is more than~1, 
then $\sfB$ is determined only up to a cyclic permutation, 
and the argument of the complex number $\alpha$ 
is different for different choices.  
That is why we give our data only for~$|\alpha|$.

%%%%%%%%%%%%%%%%%%%%%%%%%%%%%%%%%%%%%%%%%%%%%%%%%%%%%%%%%%%%
%%%%%%%%%%%%%%%%%%%%%%%%%%%%%%%%%%%%%%%%%%%%%%%%%%%%%%%%%%%%

\section{H\"older regularity and scaling properties 
        of the boundaries of Siegel disks 
        -- numerical methods and results\label{sec:num-regul}}

\subsection{Some results from harmonic analysis\label{sec:num-theor}}  

Let $\regul = n + \regul'$, where $n\in\{0,1,2,\ldots\}$, 
and $\regul'\in(0,1)$.  
We say that a function $\phi:\bbT\to\bbT$ 
has {\em (global) H\"older regularity} $\regul$ 
and write $\phi\in C^{\regul}(\bbT)$ 
if $\regul=n+\regul'$ is the largest number for which 
$\phi^{(n)}$ exists and for some constant $C>0$ 
satisfies 
\[
\left|\phi^{(n)}(t)-\phi^{(n)}\left(s\right)\right| 
\leq C\,\left|t-s\right|^{\regul'} 
\qquad \forall\ t,s \in \bbT \ .
\]
We call attention to the fact that 
we do not allow $\regul$ to be an integer 
since otherwise the definition of H\"older so 
that the characterizations we discuss later must be modified. 
In our problem, $\regul$ turns out to be non-integer, so 
that the characterizations we discuss apply. 

In the mathematical literature, there are many characterizations of 
the H\"older regularity of functions.  
Some surveys that we have found useful are \cite{Stein1970, Krantz83}. 

In the paper \cite{LlaveP02}, we developed implementations of several 
criteria for determining H\"older regularity numerically 
based on harmonic analysis, 
and assessed the reliability of these criteria. 
In this paper, we only use one of them, 
namely the method that we called the Continuous Littlewood-Paley 
(CLP) method which has been used in 
\cite{Carletti03,ApteLP05,FuchssWAM06,OlveraP06}.  
The CLP method is based on the following theorem \cite{Stein1970,Krantz83}:  

\begin{theorem}\label{thm:CLP}
A function $\phi\in C^{\regul} (\mathbb{T})$ if and only if 
for some  $\eta>0$ there exists a constant $C>0$ 
such that for all $\tau>0$ and $\eta\geq 0$ 
\begin{equation}  \label{eq:CLP}
\left\| 
\left( \frac{\partial}{\partial \tau} \right)^\eta 
\mathrm{e}^{-\tau\sqrt{-\Delta}} \, \phi 
\right\|_{L^\infty(\bbT)} 
\leq 
C \, \tau^{\regul-\eta} \ ,
\end{equation}
where $\Delta$ stands for the Laplacian: 
$\Delta\phi(t)=\phi''(t)$.  
\end{theorem}

Note that one of the consequences of Theorem~\ref{thm:CLP}
is that if the bounds \eqref{eq:CLP} hold for some $\eta$,
they hold for any other $\eta$.  
Even if from the mathematical point of view, all values of $\eta$ 
would give the same result, it is convenient from the numerical point
of view to use several to assess the reliability of the method. 
%For rather smooth functions, the decay of the right-hand side of 
%\eqref{eq:CLP} for small $\eta$ is very fast and the quality of the fit is 
%worse. Using larger $\eta$, we can get slower decays 
%that provide more reliable fits.  

%%%%%%%%%%%%%%%%%%%%%%%%%%%%%%%%%%%%%%%%%%%%%%%%%%%%%%%%%%%%

\subsection{Algorithms used\label{sec:num-algorithms}}

The algorithm we use is based on the fact that 
$\left( \frac{\partial}{\partial \tau} \right)^\eta 
\mathrm{e}^{-\tau\sqrt{-\Delta}} 
$ is a diagonal operator when acting on a Fourier representation of 
the function: if 
\[
\phi (t) = \sum_{k\in\bbN} \widehat\phi_k \, \rme^{-2\pi i k t} \ , 
\]
then 
\[
\left( \frac{\partial}{\partial \tau} \right)^\eta 
\mathrm{e}^{-\tau\sqrt{-\Delta}} \phi(t) 
= 
\sum_{k\in\bbN} 
(-2\pi |k|)^\eta \, \rme^{-2\pi \tau |k|} \, 
\widehat\phi_k \, \rme^{-2\pi i k t} \ .
\]

The Fourier transform of the function $\phi$ 
can be computed efficiently if we are given 
the values of $\phi(t)$ on a dyadic grid, 
i.e., at the points $t_m=2^{-M}m$, 
where $M$ is some natural number 
and $m=0,1,\ldots,2^M-1$.  
Unfortunately, the computation indicated in \eqref{iters} 
gives the values of the function $\chi=\Re\chi+i\,\Im\chi$ 
on the set $\{n\sigma\}_{n=0}^N$ 
of translations by the irrational number~$\sigma$. 
Therefore, we need to perform some interpolation 
to find the approximate values of $\chi$ on the dyadic grid, 
after which we Fast Fourier Transform (FFT) 
can be computed efficiently.  

Hence, the algorithm to assess the regularity is the following.  
\begin{enumerate} 
\item 
Locate the critical point~$c$ (such that $f'(c)=0$).  
\item 
Use equation \eqref{iters} to 
obtain the values of the function $\chi$ 
at the points $\{n \sigma\}_{n = 0}^N$ for some large~$N$.  
\item 
Interpolate $\Re\chi$ and $\Im\chi$ 
to find their approximate values 
on the dyadic grid $\{2^{-M}m\}_{m = 0}^{2^M -1}$.  
\item 
For a fixed value of $\eta$, 
compute the value of the left hand side of \eqref{eq:CLP} 
for several values of 
$\tau$ by using FFT (for $\phi=\Re\chi$ and separately 
for $\phi=\Im\chi$); 
do this for several values of~$\eta$.  
\item 
Fit the decay predicted by~\eqref{eq:CLP} to find the regularity~$\regul$.  
\end{enumerate}

Let us estimate the cost in time and storage of 
the algorithm above
keeping  $2^M$ values of the  function. 
Of course, locating the critical point $c$ is trivial. 
Iteration and interpolation require $O(2^M)$ operations. 
Then, each of the calculations of \eqref{eq:CLP} requires
two FFT, which is $O(2^M \ln (2^M)) = O(M\,2^M)$. 
In the computers we used (with about $1\,\mathrm{GB}$ of 
memory) the limiting factor was the storage,
but keeping several million iterates of $f$ and computing 
$2^{23}\approx 8\times 10^6$ Fourier coefficients 
of $\Re\chi$ and $\Im\chi$ was quite feasible. 
(Note that a double precision array of $2^{23}$ 
double complex numbers takes $2^{27}$ bytes 
$= 128\,{\rm MB}$ and one needs to have several copies.)  

The iterates $f^n(c)$ were computed 
by using extended precision with GMP 
-- an arbitrary-precision extension of~$C$ language~\cite{GMP}.  
This extra precision is very useful to avoid that the orbit 
scales. Note that the critical points are at the boundary of 
the domains of stability, so that they are moderately unstable. 

The extended precision is vitally important 
in computing the scaling exponents~$\alpha$.  
To obtain each value in Table~\ref{table:exponents}, 
we computed several billion iterates of the critical point of the map.  
To reduce the numerical error, we used several hundred 
digits of precision.

%%%%%%%%%%%%%%%%%%%%%%%%%%%%%%%%%%%%%%%%%%%%%%%%%%%%%%%%%%%%

\subsection{Visual explorations\label{sec:visual-explorations}}

%%%%%%%%%%%%%%%%%%%%%%%%%%%%%%%%%%%%%%%%%%%%%%%%%%%%%%%%%%%%%%
\begin{comment}
To make the pictures of the spectra, do the following:
* Go to   SIEGEL/PAPER   on my Mac laptop.
* Switch the shell to  bash.
* Run the shell script  run_make_pictures_spectra.sh 
  with command-line argument the rotation number 
  of the map, i.e.
     run_make_pictures_spectra.sh  4_inf
\end{comment}
%%%%%%%%%%%%%%%%%%%%%%%%%%%%%%%%%%%%%%%%%%%%%%%%%%%%%%%%%%%%%%

\begin{figure}
        \includegraphics[width=0.49\textwidth]{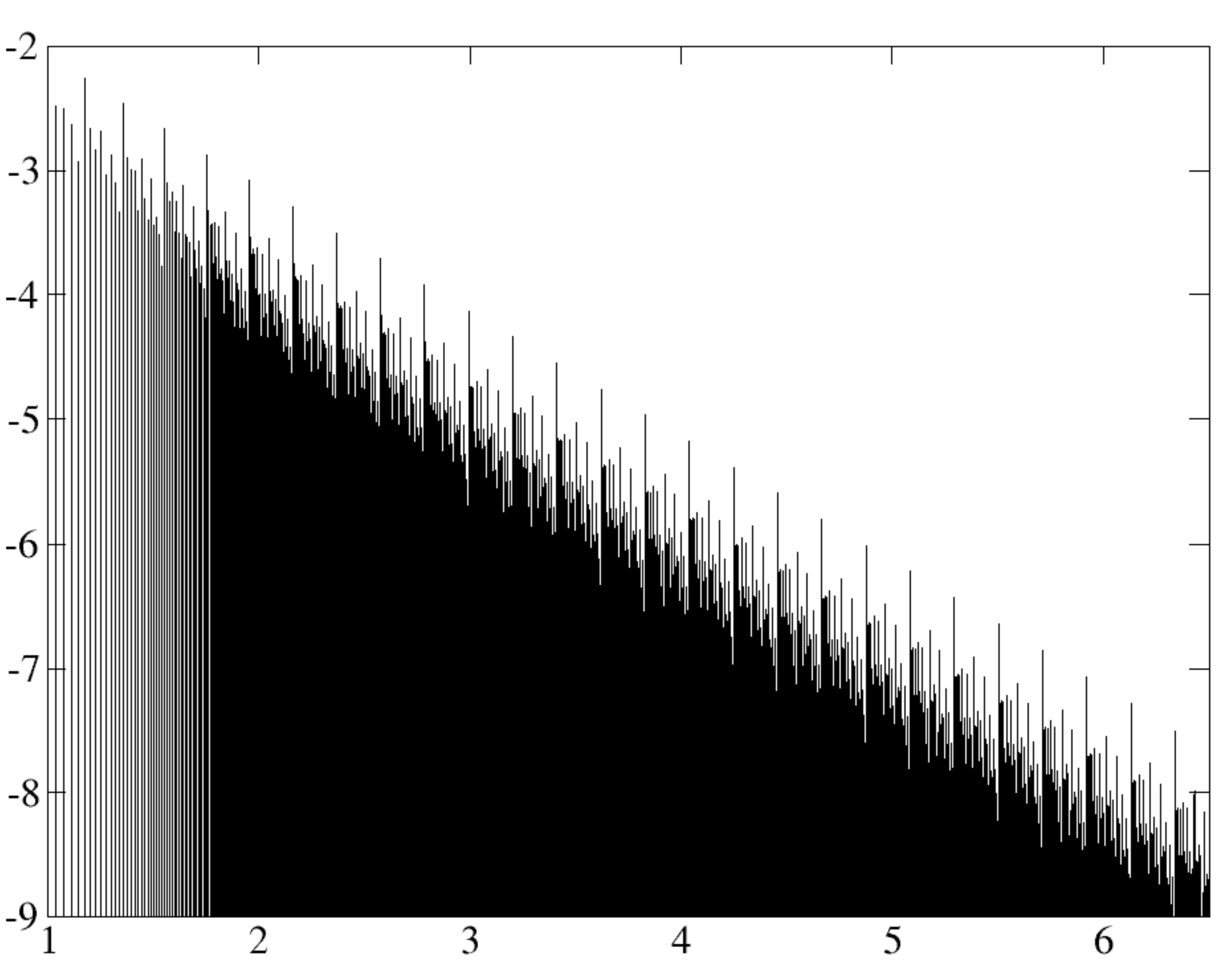}%width=0.49\textwidth
        \caption{\label{fig:spectrum_1_inf}
        Plot of $\log_{10}|\,(\widehat{\Re\chi})_k|$ vs.~$\log_{10}|k|$ 
        for the map \eqref{eq:old-f} with 
        $\sigma=\langle 1^\infty\rangle$ 
        -- see text for details.  
}
\end{figure}

\begin{figure}
        \includegraphics[width=0.49\textwidth]{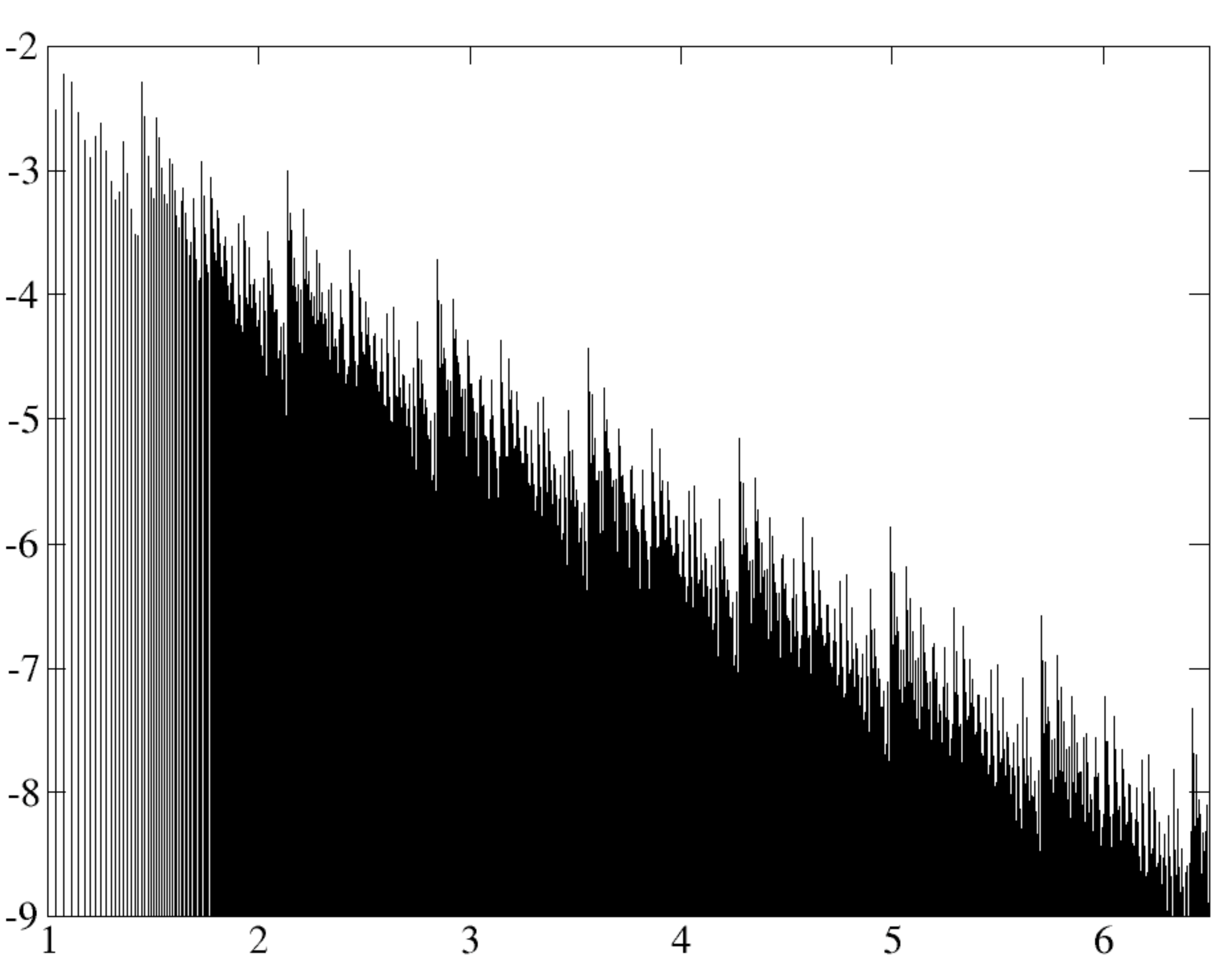}%width=0.49\textwidth
        \caption{\label{fig:spectrum_5_inf}
        Plot of $\log_{10}|\,(\widehat{\Re\chi})_k|$ vs.~$\log_{10}|k|$ 
        for the map \eqref{eq:old-f} with 
        $\sigma=\langle 5^\infty\rangle$ 
        -- see text for details.  
}
\end{figure}

\begin{figure}
        \includegraphics[width=0.49\textwidth]{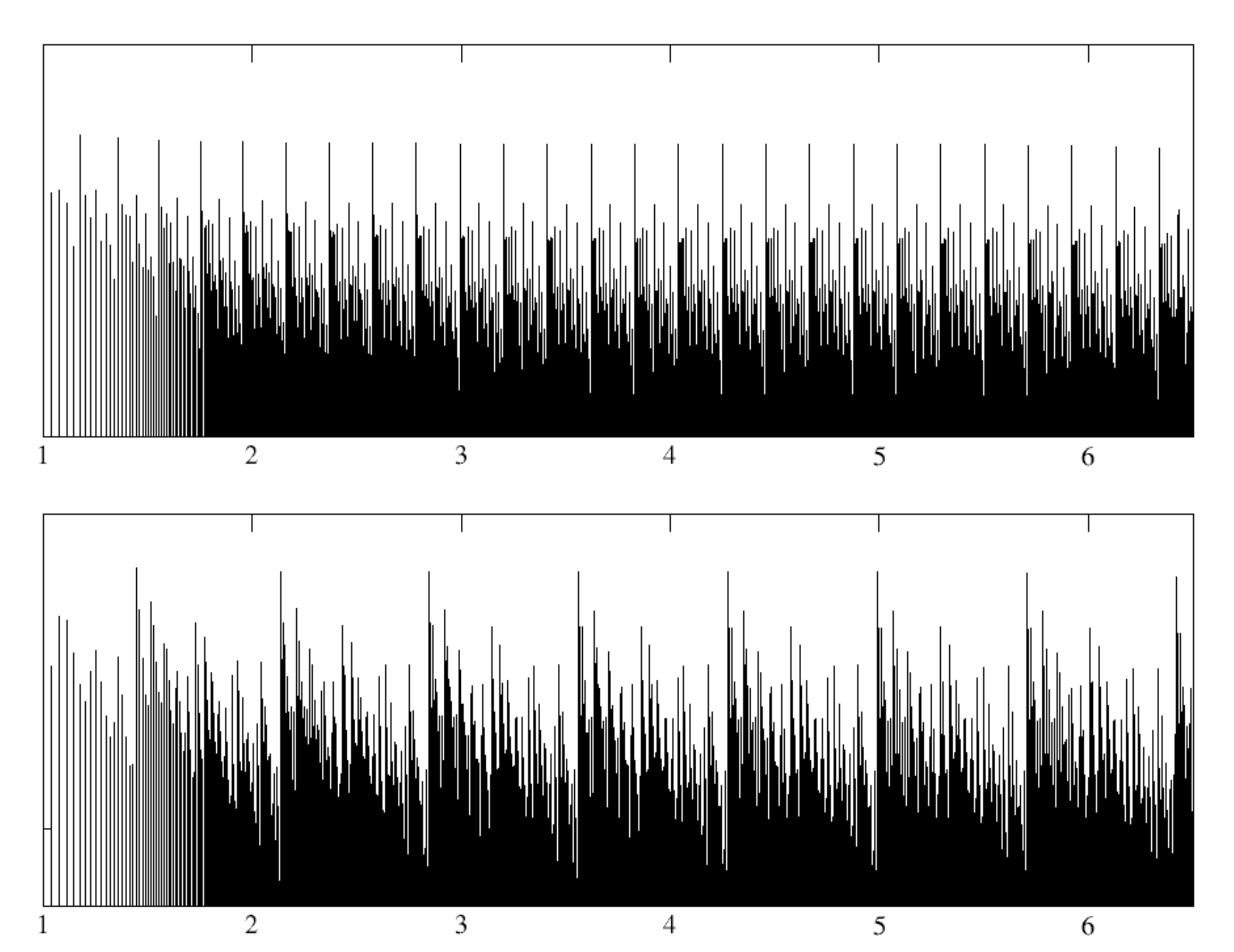}%width=0.49\textwidth
        \caption{\label{fig:spectra_1_inf_and_4_inf_straight}
        Plot of $\log_{10}|k\,(\widehat{\Re\chi})_k|$ vs.~$\log_{10}|k|$ 
        for the map \eqref{eq:old-f} with 
        $\sigma=\langle 1^\infty\rangle$ (top), 
        $\sigma=\langle 5^\infty\rangle$ (bottom) 
        -- see text for details.  
}
\end{figure}

In Figures \ref{fig:spectrum_1_inf} and \ref{fig:spectrum_5_inf} 
we have plotted (with impulses) 
the modulus of the Fourier coefficients 
$(\widehat{\Re\chi})_k$ versus $|k|$ on a log-log scale 
(for several million values of~$k$) 
for the boundary of the Siegel disk 
corresponding to the map \eqref{eq:old-f} 
for rotation numbers $\sigma$ 
equal to $\langle 1^\infty\rangle$ 
and $\langle 5^\infty\rangle$, 
and order of the critical point $d=1$ 
(the same cases as the ones in 
Figures \ref{fig:orbits_1_inf_and_5_inf_and_5_1_inf} and~\ref{fig:functions1}).  
The self-similar structure of the boundary of the Siegel 
disk is especially clearly visible 
in the ``straightened-out'' graph of the spectrum 
-- in Figure~\ref{fig:spectra_1_inf_and_4_inf_straight} 
we plotted $\log_{10}|k\,(\widehat{\Re\chi})_k|$ vs.~$\log_{10}|k|$ 
for the same spectra as in Figures 
\ref{fig:spectrum_1_inf} and~\ref{fig:spectrum_5_inf}.  
The width of each of the periodically repeating ``windows'' in the figure 
(i.e., the distance between two adjacent high peaks) 
is approximately equal to $|\log_{10}\sigma|$, 
where $\sigma$ is the corresponding rotation number.  
The periodicity in the Fourier series has been related to 
some renormalization group analysis in phase space
\cite{Shraiman84}.

\begin{figure}
        \includegraphics[width=0.5\textwidth]{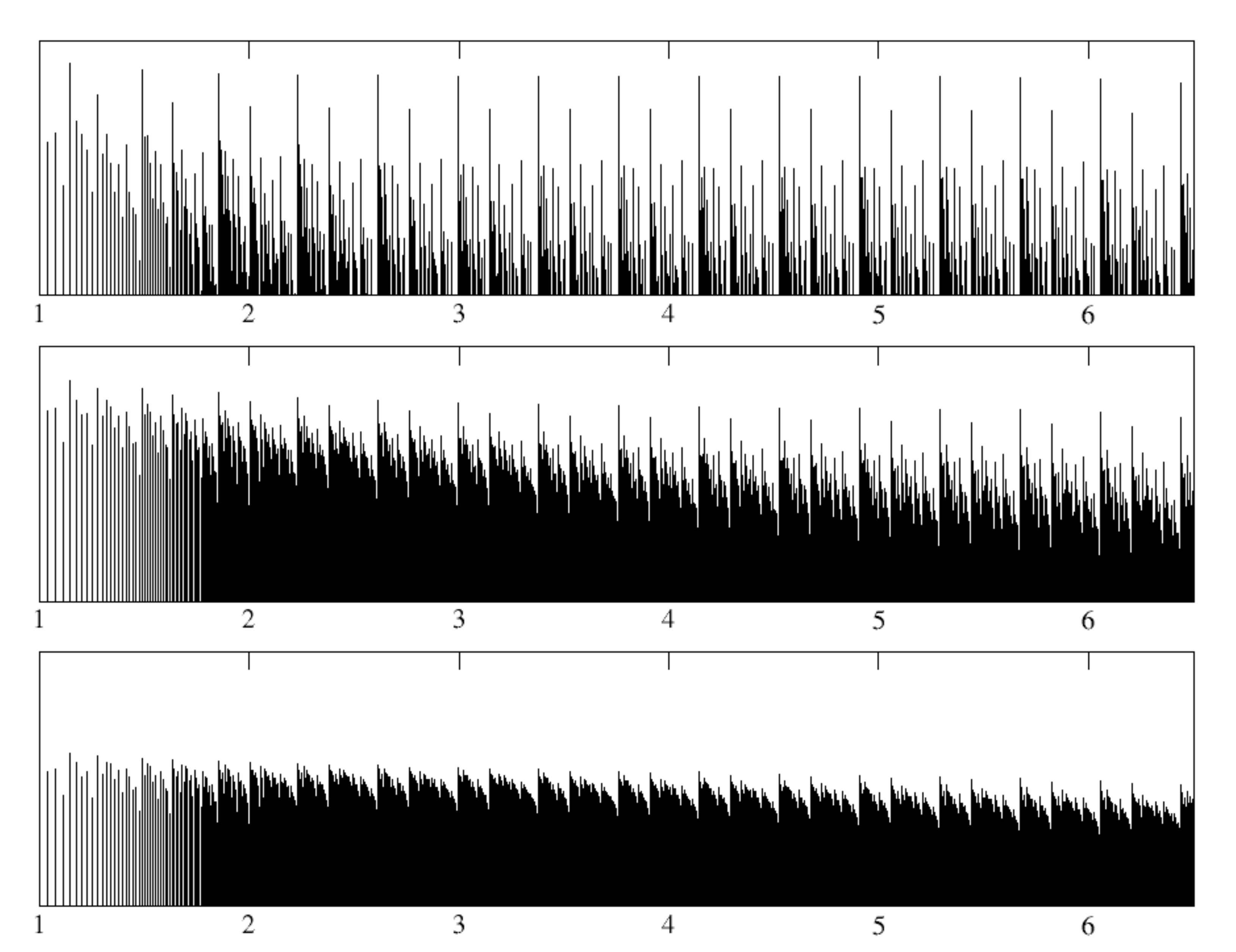}%width=0.50\textwidth
        \caption{\label{fig:spectrum_2_inf_crit_1_5_20_reduced}
        Plot of $\log_{10}|k\,(\widehat{\Re\chi})_k|$ vs.~$\log_{10}|k|$ 
        for the maps $f_{d,\langle 2^\infty \rangle, 1+3i}$ 
        with $d=1,5,20$.  
}
\end{figure}

To illustrate the effect of the order of criticality 
on the Fourier spectrum of $\Re\chi$, 
we showed in Figure~\ref{fig:spectrum_2_inf_crit_1_5_20_reduced} 
the ``straightened-out'' graphs of the spectra, 
$\log_{10}|k\,(\widehat{\Re\chi})_k|$ versus $\log_{10}|k|$, 
of the function $\Re\chi$ corresponding to the maps 
$f_{d,\langle 2^\infty \rangle, 1+3i}$ 
for orders $d=1,5,20$ 
(the Siegel disks of these maps were shown 
in Figure~\ref{fig:high-crit}).  
For all plots in this figure we used the same scale in vertical direction.  
An interesting observation 
-- for which we have no conceptual explanation at the moment 
-- is that the variability of the magnitudes of the Fourier coefficients 
decreases as the order of the critical point increases.  

\begin{figure}
%\centering
        \includegraphics[width=0.47\textwidth]{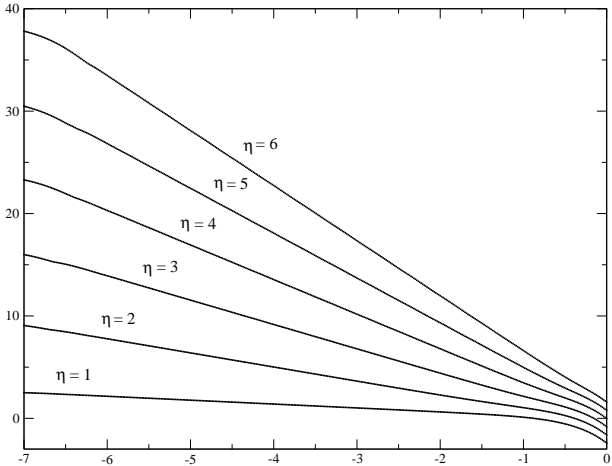}\\[2mm]%width=0.47\textwidth
        \includegraphics[width=0.48\textwidth]{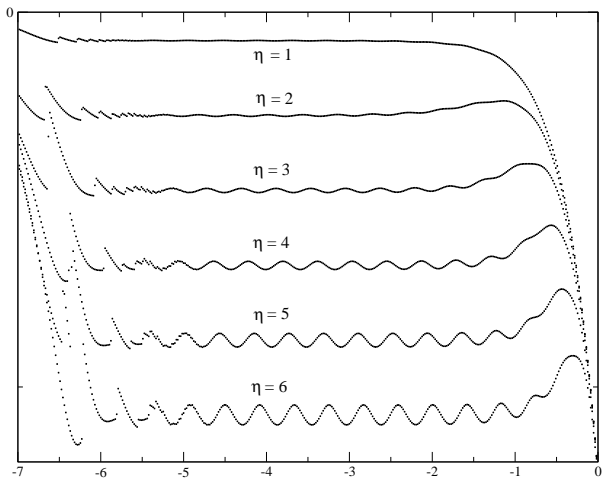}%width=0.48\textwidth
        \caption{\label{plot_CLP_x_1_inf-1}
        {\em Top:} 
        Plot of 
        $\log_{10} N_\eta(\tau)$ vs.~$\log_{10}\tau$ 
        for the map \eqref{eq:old-f} 
        with $\sigma=\langle 1^\infty\rangle$.  
        {\em Bottom:} 
        Plot of the first differences 
        between the vertical coordinates of adjacent points 
        in the top figure.  
}
\end{figure}

\begin{figure}
%%\begin{flushright}
        \includegraphics[width=0.47\textwidth]{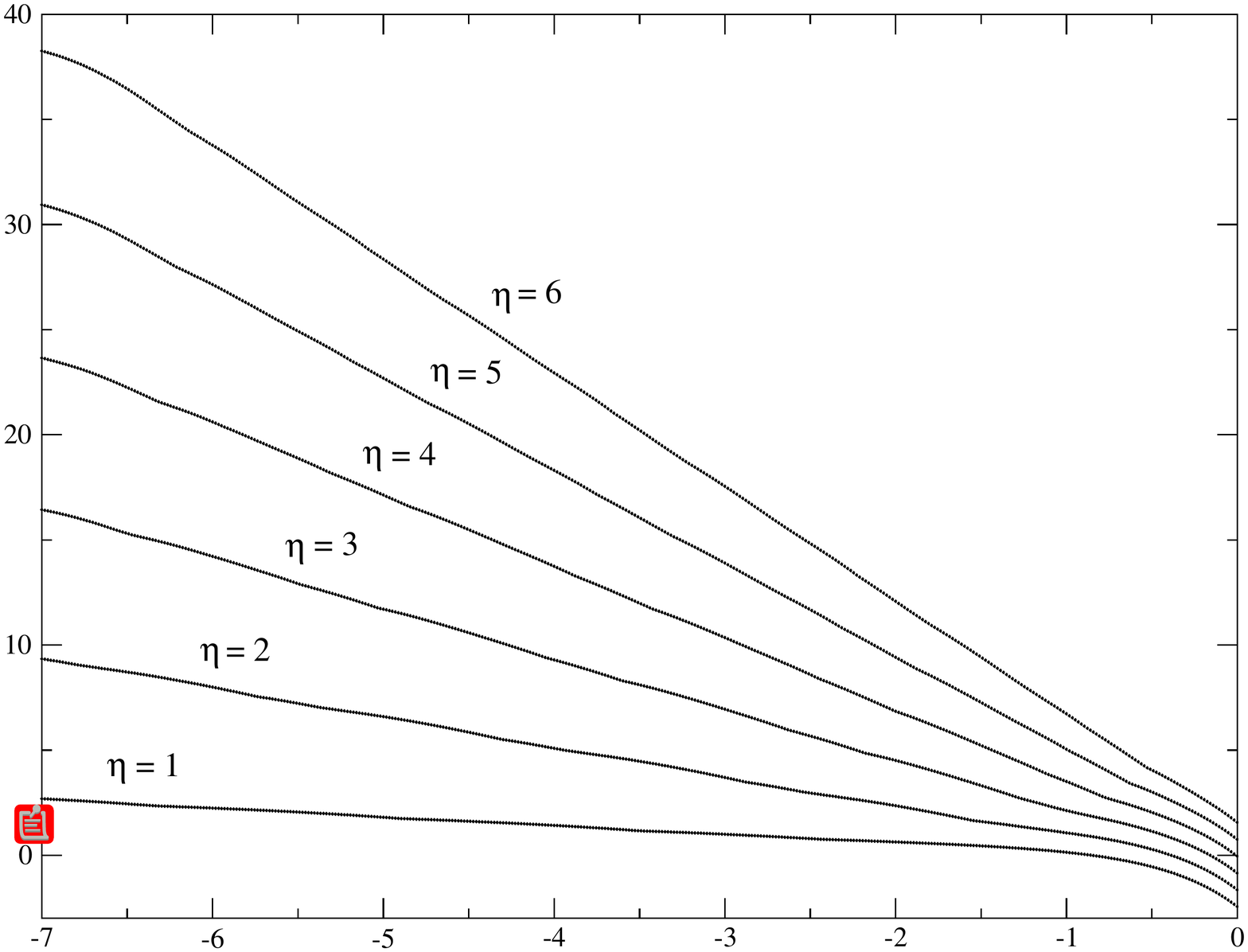}\\[2mm]%width=0.47\textwidth
        \includegraphics[width=0.48\textwidth]{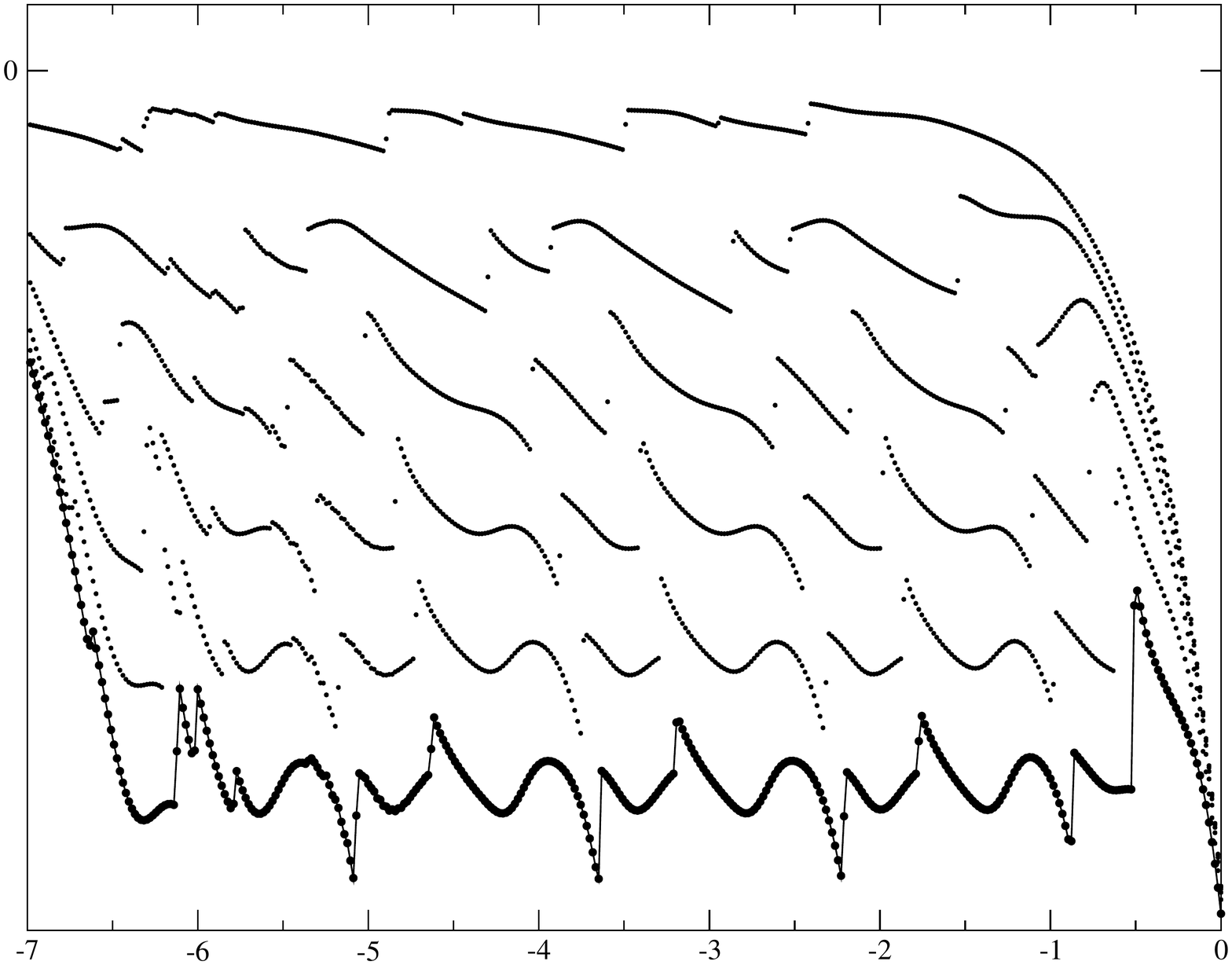}%width=0.48\textwidth
        \caption{\label{plot_CLP_x_5_inf-1}
        {\em Top:} 
        Plot of 
        $\log_{10} N_\eta(\tau)$ vs.~$\log_{10}\tau$ 
        for the map \eqref{eq:old-f} 
        with $\sigma=\langle 5^\infty\rangle$.  
        {\em Bottom:} 
        Plot of the first differences 
        between the vertical coordinates of adjacent points 
        in the top figure.  
}
%\end{flushright}
\end{figure}

Figures~\ref{plot_CLP_x_1_inf-1} and \ref{plot_CLP_x_5_inf-1} 
illustrate the CLP method (Theorem~\ref{thm:CLP}) in practice.  
In the top part of each figure we have plotted 
on a log-log scale the norms in the left-hand side 
of~\eqref{eq:CLP}, 
\begin{equation}  \label{eq:N-eta}
N_\eta(\tau):= \left\|\left(\frac{\partial}{\partial\tau}\right)^\eta 
        \, \rme^{-\tau\sqrt{-\Delta}} \, 
        \Re\chi\right\|_{L^\infty(\bbT)} \ ,
\end{equation}
as functions of $\tau$, for $\eta=1,2,\ldots,6$ 
for the map \eqref{eq:old-f} 
with for rotation numbers 
$\langle 1^\infty\rangle$ and~$\langle 5^\infty\rangle$, respectively.  
For each value of $\eta$, the ``line'' consists 
of 400 points corresponding to 400 different values 
of $\tau$ for which we have computed the corresponding norm.  

The bottom parts of 
Figures~\ref{plot_CLP_x_1_inf-1} and \ref{plot_CLP_x_5_inf-1} 
show the behavior of the first differences 
between the vertical coordinates of adjacent points 
from the top parts of the figures.  
Clearly, the points in the top parts 
do not lie on exact straight lines 
but have small periodic (as functions of $\log \tau $) 
displacements.  
To make this more clear in the bottom part of 
Figure~\ref{plot_CLP_x_5_inf-1}, 
we have plotted the points for $\eta=6$ 
with larger circles, 
and have connected them with lines.  
%The plot of $\eta=6$ for rotation number 
%$\langle 5^\infty\rangle$ 
%is given with larger circles and lines 
%connecting them to show the ``periodic'' behavior 
%of the first differences.  

We would like to point out that the bound \eqref{eq:CLP} 
%does not guarantee that, for a fixed $\eta$, 
%the norm in the left-hand side of \eqref{eq:CLP} 
%should decay approximately as $\tau^{\regul-\eta}$ 
%because \eqref{eq:CLP} 
is only an upper bound.  
However, the discrete scaling invariance 
of the Siegel disk at small scales 
-- manifested also by  the existence of the ``periodic windows'' 
in Figure~\ref{fig:spectra_1_inf_and_4_inf_straight} 
-- implies that the leading behavior of 
$\log_{10} N_\eta(\tau)$ as a function of $\log_{10}\tau$ is indeed linear, 
and that \eqref{eq:CLP} is close to being saturated.  

The discrete scaling invariance at small scales allows for 
$\log_{10} N_\eta(\tau)$ as a function of $\log_{10}\tau$ 
to have small periodic corrections 
superimposed on the leading linear behavior.  
Our calculations are precise enough that these corrections 
to the leading behavior are clearly visible. 
The bottom parts of Figures \ref{plot_CLP_x_1_inf-1} 
and \ref{plot_CLP_x_5_inf-1} 
show exactly these corrections.  

The form of these small periodic corrections 
depends in a complicated way on the the behavior of the Fourier coefficients 
in the ``periodic windows'' in the Fourier spectrum.  
The presence of this periodic corrections is a good indicator 
of the ranges of $t$ which are large enough that the asymptotic
behaviour has started to take hold, but small enough so that they 
are not dominated by the round-off and truncation error. 
In our previous works \cite{LlaveP02,ApteLP05}, 
we have also found periodic corrections to the 
scaling in other conjugacies related to the 
regularity of conjugacies of other critical objects.

%%%%%%%%%%%%%%%%%%%%%%%%%%%%%%%%%%%%%%%%%%%%%%%%%%%%%%%%%%%%

\subsection{Numerical values of the H\"older regularity\label{sec:regul}}  

In Table~\ref{table:regularity} we give the computed values 
%%%%%%%%%%%%%%%%%%%%%%%%%%%%%%%%%%%%%%%%%%%%%%%%%%%%%%%%%%%%
\begin{comment}
\begin{table}
\caption{THESE VALUES ARE FOR THE FILES ``ocrit''!!!!!
}
\begin{ruledtabular}
\begin{tabular}{rccccc}
$d$ & $\langle 1^\infty\rangle$ & $\langle 2^\infty\rangle$ & $\langle 3^\infty\rangle$ & $\langle 4^\infty\rangle$ & $\langle 5^\infty\rangle$\\
\hline
 1 & 0.621 & 0.616 & 0.607 & 0.594 & 0.579 \\
 2 & 0.433 & 0.428 & 0.418 & 0.404 & 0.388 \\
 3 & 0.330 & 0.325 & 0.316 & 0.305 & 0.291 \\
 4 & 0.266 & 0.262 & 0.254 & 0.245 & 0.233 \\
 5 & 0.223 & 0.219 & 0.212 & 0.204 & 0.194 \\
 6 & 0.191 & 0.188 & 0.182 & 0.175 & 0.166 \\
10 & 0.123 & 0.121 & 0.116 & 0.111 & 0.105 \\
15 & 0.085 & 0.083 & 0.080 & 0.077 & 0.074 \\
20 & 0.065 & 0.064 & 0.061 & 0.059 & 0.054 \\
\end{tabular}
\end{ruledtabular}
\end{table}
\end{comment}
%%%%%%%%%%%%%%%%%%%%%%%%%%%%%%%%%%%%%%%%%%%%%%%%%%%%%%%%%%%%
\begin{table}
\caption{Global H\"older regularity $\regul$ for 
        maps with rotation numbers 
        $\sigma=\langle k^\infty\rangle$ ($k=1,2,3,4,5$) 
        and critical point of order~$d$.  
        \label{table:regularity}
}
\begin{ruledtabular}
\begin{tabular}{rccccc}
$d$ & $\langle 1^\infty\rangle$ & $\langle 2^\infty\rangle$ & $\langle 3^\infty\rangle$ & $\langle 4^\infty\rangle$ & $\langle 5^\infty\rangle$\\
\hline
 1 & 0.621 & 0.617 & 0.607 & 0.596 & 0.578 \\
 2 & 0.432 & 0.427 & 0.417 & 0.404 & 0.388 \\
 3 & 0.328 & 0.324 & 0.313 & 0.300 & 0.291 \\
 4 & 0.263 & 0.260 & 0.252 & 0.244 & 0.232 \\
 5 & 0.220 & 0.217 & 0.210 & 0.203 & 0.193 \\
 6 & 0.189 & 0.186 & 0.180 & 0.174 & 0.163 \\
10 & 0.121 & 0.120 & 0.115 & 0.111 & 0.105 \\
15 & 0.084 & 0.082 & 0.079 & 0.077 & 0.074 \\
20 & 0.064 & 0.063 & 0.061 & 0.058 & 0.055 \\
\end{tabular}
\end{ruledtabular}
\end{table}
of the global H\"older regularities 
of the real and imaginary parts of the dynamically natural 
parameterizations $\chi$ 
\eqref{particular}, \eqref{chi0} 
of the boundaries of the Siegel disks.  
We studied maps of the form \eqref{eq:f-def}, 
with different values of $\beta$, 
and with different orders $d$ of the critical point 
in~$\partial\calD$; 
some runs with the map \eqref{eq:old-f} 
(for which the critical point is simple) 
were also performed.  

To obtain each value in the table, we performed the procedure 
outlined in Section~\ref{sec:num-algorithms} 
for at least two maps of the form~\eqref{eq:f-def}.  
For each map we plotted the points from the CLP analysis 
for $\eta=1,2,\ldots,6$ 
as in  the top parts of Figures \ref{plot_CLP_x_1_inf-1} 
and \ref{plot_CLP_x_5_inf-1}, looked at the differences 
between the vertical coordinates of adjacent points 
(i.e., at graphs like in the bottom parts 
of Figures \ref{plot_CLP_x_1_inf-1} 
and \ref{plot_CLP_x_5_inf-1}), 
and selected a range of values of $\log_{10}\tau$ 
for which the differences oscillate regularly.  
For this range of $\log_{10}\tau$, 
we found the rate of decay of the norms 
$N_\eta(\tau)$ \eqref{eq:N-eta} 
by measuring the slopes, $\regul-\eta$, 
from which we computed the regularity~$\regul$.  

The accuracy of these values is difficult to estimate, 
but a conservative estimate on the relative error 
of the data in Table~\ref{table:regularity} 
is about~3\%.  

We have also computed the regularity $\regul$ 
and the scaling exponent $\alpha$ of several maps 
with rotation numbers with the same tails  
of the continued fraction expansion
but with different heads.  

Within the accuracy of our computations, 
the results did not depend on the head, 
which is consistent with the predictions 
of the renormalization group picture.  
%we give in Table~\ref{table:regularity} 
%the computed values of the regularities 
%for several rotation numbers 
%with the same tails but different heads.  
%The accuracy of our calculations  did not allow us to resolve 
%the regularities for $\sigma=\langle (1^52)^\infty \rangle$, 
%$\langle (1^62)^\infty \rangle$, $\langle (1^72)^\infty \rangle$ 
%and $\langle (1^82)^\infty \rangle$ (and all their cyclic permutations).  

%%%%%%%%%%%%%%%%%%%%%%%%%%%%%%%%%%%%%%%%%%%%%%%%%%%%%%%%%%%%

\subsection{Importance of the phases of Fourier coefficients\label{sec:phases}}

In Figure~\ref{fig:spectra_1_inf_and_4_inf_straight} 
we saw that the modulus of the Fourier coefficients 
of $\Re\chi$ and $\Im\chi$ decreases very approximately as 
\begin{equation}  \label{eq:decay-modulus}
\left| (\widehat{\Re\chi})_k\right| \leq \frac{C}{|k|} \ , 
\qquad 
\left| (\widehat{\Im\chi})_k\right| \leq \frac{C}{|k|} \ .
\end{equation}
and that these bounds come close to be saturated infinitely 
often. This seems to be true independent of what 
the rotation number is.

This rate of decay of the Fourier coefficients 
would be implied by $\Re\chi$ and $\Im\chi$ being $C^1$, 
but from \eqref{eq:decay-modulus} 
we cannot conclude that $\Re\chi$ and $\Im\chi$ 
are even continuous 
(note, for example, 
that the function $f(x)=\sum_{k=1}^\infty \frac1k \cos kx$ 
is discontinuous at $x=0$).  
It is well-known from harmonic analysis that the phases 
of the Fourier coefficients 
play a very important role, 
and changing the phases of Fourier coefficients 
changes the regularity of the functions (see, 
e.g., \cite{Zygmund02}) 
\[
\sum_{k=2}^\infty \, \frac1{k^{1/4}} \,\, e^{2\pi i kx} 
\sim
\frac1{|x|^{3/4}}  \quad \mbox{as} \ |x|\to 0  \ ,
\]
while 
\[
\sum_{k=2}^\infty \, \frac{e^{i\sqrt{k}}}{k^{1/4}} \,\, e^{2\pi  kx} 
\sim 
\frac1{|x|^{2}}  \quad \mbox{as} \ |x|\to 0 \ .
\]

\begin{figure}
        \includegraphics[width=0.48\textwidth]{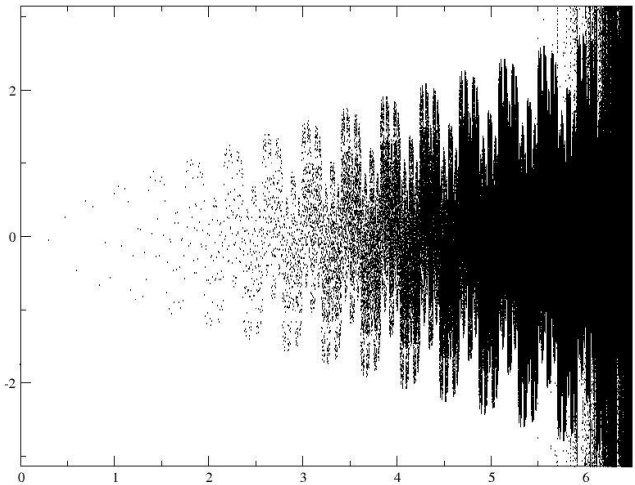}%width=0.48\textwidth
        \caption{\label{fig:phases}
        Plot of the phases of the Fourier coefficients 
        (vs.~$\log_{10}|k|$) 
        of $\Re\chi$ for the map~$f_{0,\langle 1^\infty \rangle,1}$.  
}
\end{figure}

In Figure~\ref{fig:phases} we depict 
the phases of the Fourier coefficients of $\Re\chi$ 
for the map $f_{0,\langle 1^\infty \rangle,1}$ 
%(in the notations of~\eqref{eq:f-def}), 
whose only critical point, $c=1$, 
is simple (i.e., of order $d=1$).  
We see that, for small $k$, the phases have a repeated pattern.  
If we consider $k\in I_j := [\sigma^{-j},\sigma^{-j-1}]$ 
(where $\sigma=\langle 1^\infty \rangle$ is the golden mean), 
we see that the phase restricted to $I_j$ 
has a pattern very similar to $I_{j+1}$, 
except that the latter is reversed and amplified.  
Of course, since the phase only takes values between 
$-\pi$ and $\pi$, the amplification of the patterns 
can only be carried out a finite number of times 
until the absolute values of the phases reach $\pi$, 
after which they will start ``wrapping around'' 
the interval $[-\pi,\pi]$.  
Unfortunately, to see this effect numerically, 
we would need hundreds of millions of Fourier coefficients, 
which at the moment is out of reach.

\begin{figure}[h]
        \includegraphics[width=0.48\textwidth]{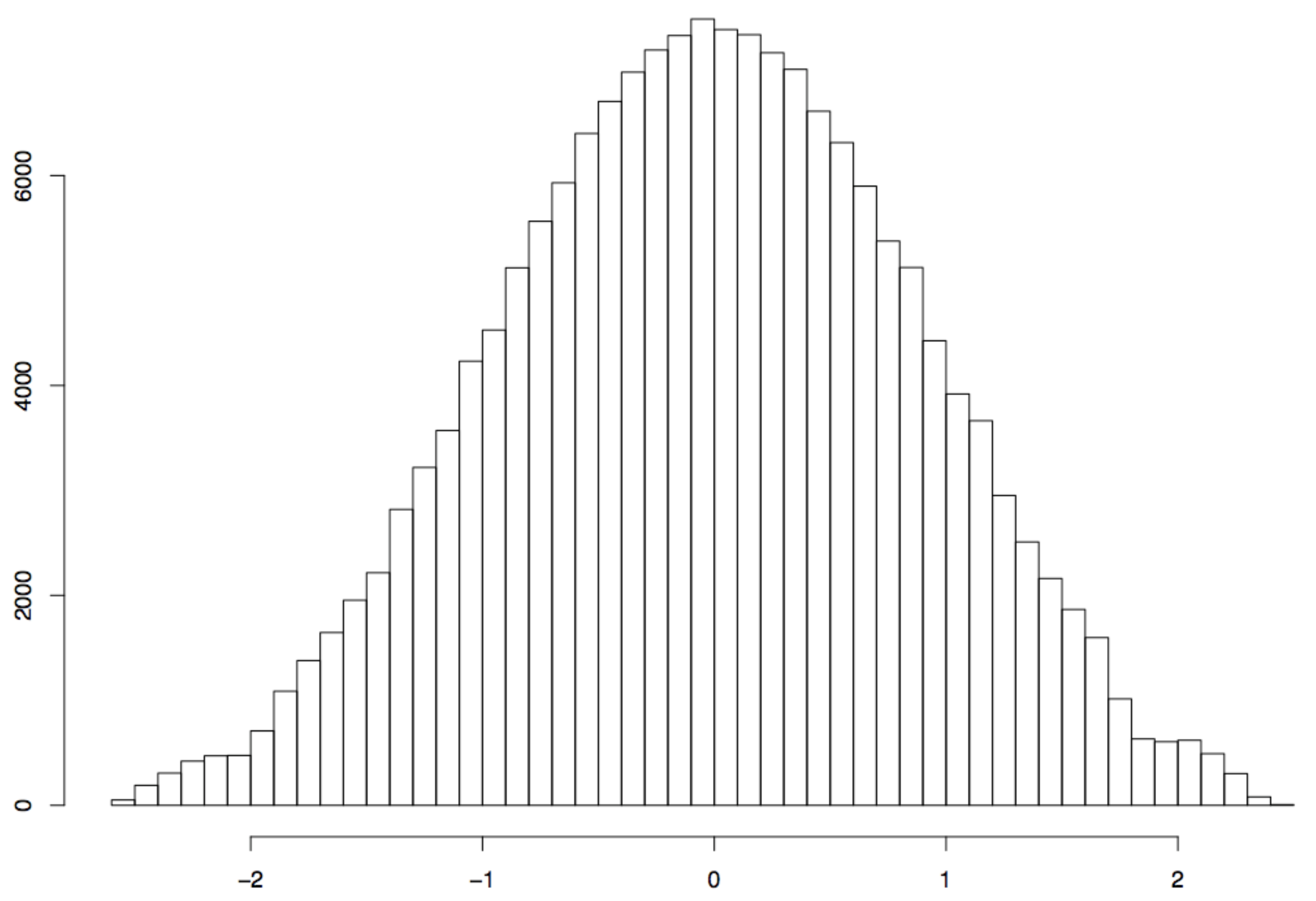}%width=0.48\textwidth
        \caption{\label{fig:phases-histo}
                Histogram of some of the phases shown 
        in Figure~\ref{fig:phases} (from $I_{24}\cup I_{25}$).  
        }
\end{figure}

\begin{figure}[h]
        \includegraphics[width=0.48\textwidth]{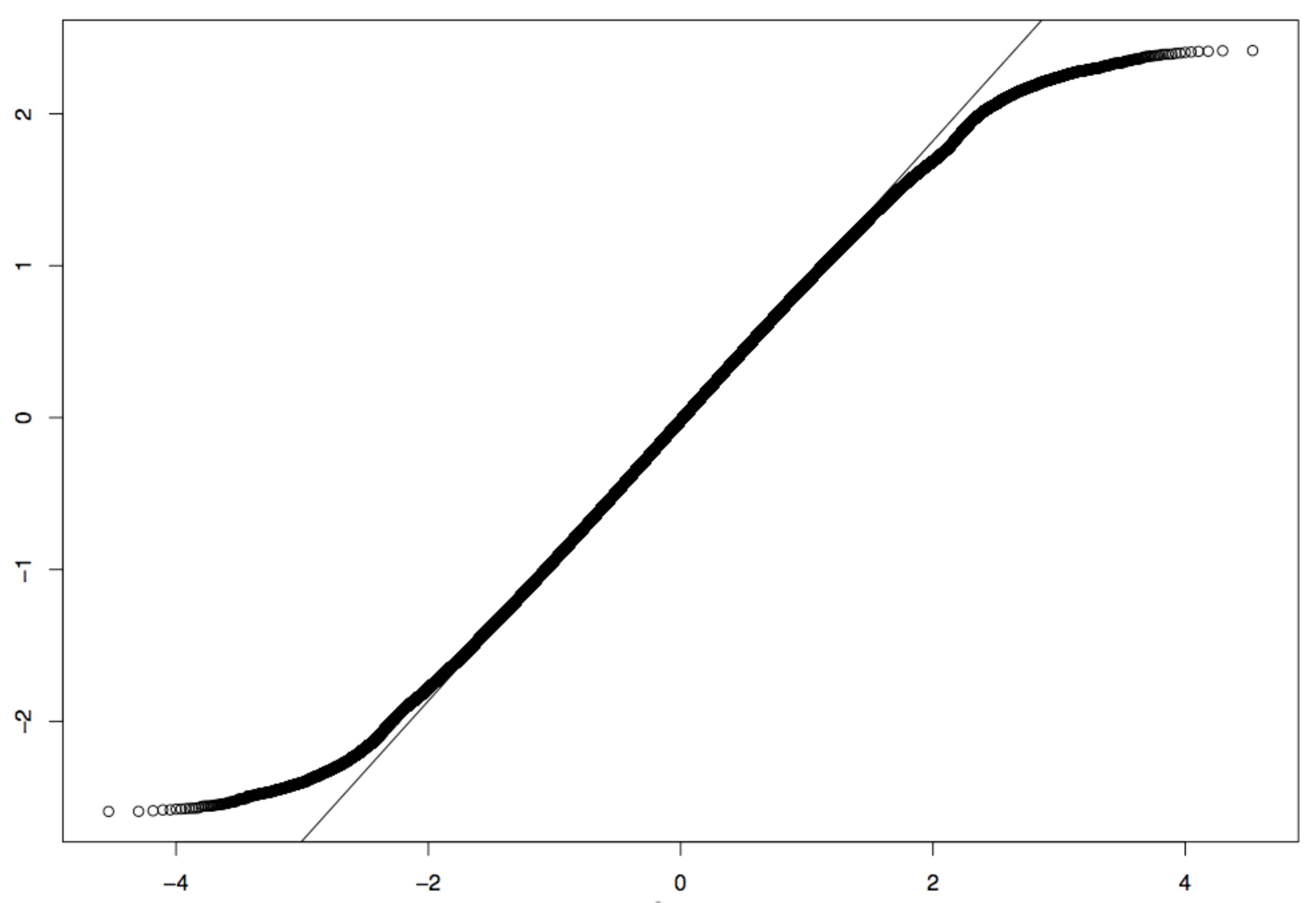}%width=0.48\textwidth
        \caption{\label{fig:phases-ks}
        Kolmogorov-Smirnov test for normality 
        of the phases in $I_{24}\cup I_{25}$.  
        }
\end{figure}

Given the above observation, 
it is natural to study the distribution 
of the phases of the Fourier coefficients  in an interval of
self-similarity. 
In Figure~\ref{fig:phases-histo} we present the histogram 
of the phases in the interval 
$k\in I_{24}\cup I_{25}=[\sigma^{-24},\sigma^{-26}]$ 
(i.e., of about 170,000 phases).  
We note that the histogram is very similar to a Gaussian.  
This visual impression is confirmed 
by using the Kolmogorov-Smirnov test, 
shown in Figure~\ref{fig:phases-ks}.  
Recall that the Kolmogorov-Smirnov test 
consists in plotting the empirical distribution 
versus the theoretical one 
(for details see, e.g., \cite[Chapter~7]{Sheskin07}).  
If indeed the empirical distribution was a sample 
of the theoretical distribution, 
we would get a set of points close the diagonal.  
The Kolmogorov-Smirnov test is available 
in many statistical packages 
(we used the package R, in which 
the command  \verb!qqnorm! gives a KS-test and the 
command \verb!qqline! displays the result of a fitted Gaussian).  
The Kolmogorov-Smirnov test reveals that, 
as expected (since the variable is an angle), 
the distribution of the phases 
has discrepancies with a  Gaussian near the edges, 
$-\pi$ and~$\pi$. Nevertheless, there is a remarkably 
good fit away from these edges. For the intervals 
we chose, most of the data points are indeed out of 
the edges.

\begin{comment}
Here are the lengths of the files with the phases:

\begin{verbatim}
   3300    3300   31350 data2
   7812    7812   74139 data3
  21491   21491  204139 data4
  67656   67656  643016 data5
 173093  173093 1645841 data6
 400659  400659 3808066 data7
\end{verbatim}

\end{comment}

%%%%%%%%%%%%%%%%%%%%%%%%%%%%%%%%%%%%%%%%%%%%%%%%%%%%%%%%%%%%

\subsection{Data on the scaling exponents\label{sec:data-exponents}}

In Table~\ref{table:exponents} 
\begin{table*}
\caption{Values of $|\alpha|$ for maps of the form $f_{d,\sigma,\beta}$ 
        \eqref{eq:f-def} 
        with different rotation numbers $\sigma$ and orders $d$ 
        of the critical point.\label{table:exponents}
}
\begin{ruledtabular}
\begin{tabular}{rlllll}
\mbox{$d$} & \mbox{$\langle 1^\infty\rangle$}  & \mbox{$\langle 2^\infty\rangle$} & \mbox{$\langle 3^\infty\rangle$} & \mbox{$\langle 4^\infty\rangle$} & \mbox{$\langle 5^\infty\rangle$} \\
\hline
 1  & 0.74193223170 & 0.5811130545 & 0.484541021 & 0.424632459 & 0.385769294 \\
 2  & 0.81215810740 & 0.686013947  & 0.607281233 & 0.55822367  & 0.5268809   \\
 3  & 0.853450202   & 0.7508249    & 0.6852424   & 0.6441419   & 0.6179964   \\
 4  & 0.88014575    & 0.793968     & 0.738015    & 0.7027580   & 0.680345    \\
 5  & 0.8987131     & 0.824557     & 0.775859    & 0.7450319   & 0.725425    \\
 6  & 0.912340      & 0.847314     & 0.804246    & 0.7768796   & 0.759458    \\
10  & 0.943087      & 0.89962      & 0.87026     & 0.851416    & 0.839375    \\
15  & 0.960463      & 0.92977      & 0.90881     & 0.895266    & 0.88658     \\
20  & 0.969717      & 0.94601      & 0.92971     & 0.919149    & 0.91236     \\
40  & 0.984450      & 0.97211      & 0.96356     & 0.9579      & 0.9544      \\
60  & 0.989504      & 0.981138     & 0.975314    & 0.97151     & 0.96906     \\
80  & 0.9920788     & 0.98575      & 0.981334    & 0.97845     & 0.97658     \\
100 & 0.993638      & 0.98854      & 0.98499     & 0.98267     & 0.98116     \\
200 & 0.996793      & 0.99422      & 0.99242     & 0.991241    & 0.99048     \\
300 & 0.997857      & 0.99613      & 0.99493     & 0.994140    & 0.99363     
\end{tabular}
\end{ruledtabular}
\end{table*}
we give the values of the modulus of the scaling exponent $\alpha$ 
for maps of the form \eqref{eq:old-f} with orders 
$d=1,\ldots,6,10,15,20,40,60,80,100$ 
of the critical point and rotation numbers $\sigma=\langle k^\infty \rangle$ 
with $k=1,\ldots,5$.  
We believe that the numerical error in these values 
does not exceed 2 in the last digit.  

\begin{figure}[h]
        \includegraphics[width=0.49\textwidth]{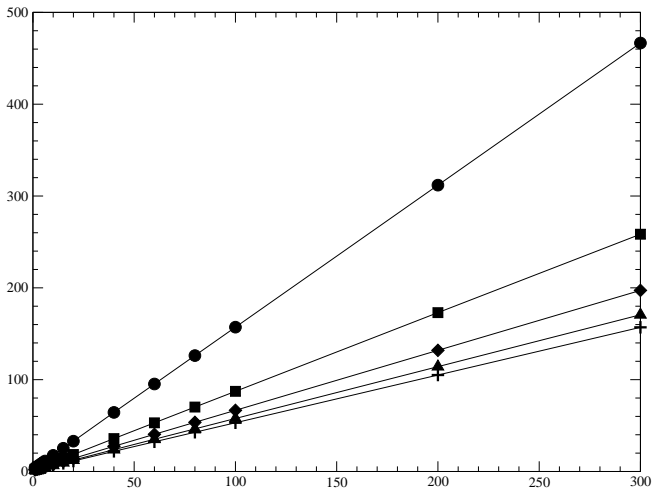}%width=0.47\textwidth
        \caption{Plot of $1/(1-|\alpha|)$ vs.~$d$ 
        for rotation numbers 
        $\langle 1^\infty \rangle$ (circles), 
        $\langle 2^\infty \rangle$ (squares), 
        $\langle 3^\infty \rangle$ (diamonds), 
        $\langle 4^\infty \rangle$ (triangles), 
        $\langle 5^\infty \rangle$ (pluses).  
        \label{fig:slopes_exponents}
        }
\end{figure}

In \cite{Osbaldestin92}, the author computed 
the scaling exponents for maps 
with rotation number~$\langle 1^\infty \rangle$ 
and critical point of different orders $d$, 
and suggested that the behavior of $\alpha$ 
for large $d$ is 
\begin{equation}  \label{eq:A}
|\alpha_{\langle 1^\infty \rangle,d}| 
\sim 
1-\frac{A_{\langle 1^\infty \rangle}}{d} 
\qquad \mbox{ as } \ d\to\infty \ .
\end{equation}
We studied the same problem for other rotation numbers
and, taking advantage
of the extended precision,  we carried out 
the computation for rather  
high degrees of criticality ($\approx 300$), see 
Table~\ref{table:exponents}.  
In Figure~\ref{fig:slopes_exponents} 
we plotted $1/(1-|\alpha|)$ versus $d$ 
for five rotation numbers.  
Our data that for high values of $d$ 
the modulus of $\alpha$ tends to~1 
for any rotation number.  
The values of the constants $A_\sigma$ 
in \eqref{eq:A} for rotation numbers 
$\langle k^\infty \rangle$ with $k=1,2,3,4,5$ 
are approximately $0.646$, $1.168$, $1.531$, $1.960$, $1.925$, 
respectively 
(the linear regression was based on the values 
for $d=40,\ldots,300$).  
%%%%%%%%%%%%%%%%%%%%%%%%%%%%%%%%%%%%%%%%%%%%%%%%%%%%%%%%%%%%
\begin{comment}
Here is the data used to plot this figure:
  1    3.8750    2.3873    1.9400    1.7380    1.6281
  2    5.3236    3.1849    2.5464    2.2636    2.1136
  3    6.8236    4.0132    3.1770    2.8101    2.6178
  4    8.3435    4.8536    3.8170    3.3643    3.1284
  5    9.8729    5.6999    4.4615    3.9221    3.6420
  6   11.4077    6.5494    5.1085    4.4819    4.1573
 10   17.5707    9.9621    7.7077    6.7302    6.2257
 15   25.2908   14.2389   10.9661    9.5480    8.8168
 20   33.0218   18.5219   14.2268   12.3684   11.4103
 40   64.3087   35.8551   27.4424   23.7530   21.9298
 60   95.2744   53.0166   40.5088   35.1000   32.3206
 80  126.2435   70.1754   53.5733   46.4037   42.6985
100  157.1833   87.2600   66.6223   57.7034   53.0786
200  311.8179  173.0104  131.9261  114.1683  105.0420
300  466.6356  258.3979  197.2387  170.6485  156.9859
\end{comment}
%%%%%%%%%%%%%%%%%%%%%%%%%%%%%%%%%%%%%%%%%%%%%%%%%%%%%%%%%%%%

\section{Calculation of the Siegel radius and the area 
of the Siegel disk\label{sec:radius-area}}

As a byproduct of our calculations we can obtain 
rather precise values of two quantities of 
mathematical interest: the area of the Siegel disk and 
the Siegel radius. 

\subsection{Calculation of Siegel radius\label{sec:radius}} 

We note that the parameterization, $\chi$, of the boundary $\partial\siegeldisk$ is 
related to the conjugacy $h$ \eqref{conjugacy} by \eqref{particular}.  
%just $\chi(s) = h(r_s \exp( i s) )$ where $h$ is 
%the conjugacy satisfying \eqref{conjugacy}. 
Then, the Fourier coefficients of $\chi$ 
satisfy $|\hat \chi_{-k}| = |h_k| \, r_\mathrm{S}^k$ 
for $k\in\bbN$, 
where $h_k$ are the Taylor coefficients of $h$ 
(recall that $h_0=0$ and $h_1=1$).  
As shown in \cite{SiegelM95, Llave01c}, one can 
get the all the coefficients $h_k$ by equating 
terms of like powers in \eqref{conjugacy}, 
and this gives infinitely many different ways to compute 
$r_\mathrm{S}$. 
In particular, since  $h_1 = 1$, 
we have $r_\mathrm{S} = |\hat \chi_{-1}|$. Since we also 
have $h_2 = f_2/[a(a -1)]$ and $|a|=1$, we obtain 
$r_\mathrm{S}^2 = |a-1| |\hat\chi_{-2}| / |f_2|$ 
(where $f_1=a$, $f_2$, $\ldots$ are the Taylor coefficients 
of the function~$f$).  
%$r_\mathrm{S}^2 = f_2/(a^2 -1) \hat \chi_k$.
Similar formulas for higher order terms are also 
available. 

\subsection{Calculation of the area of the 
Siegel disk} 

Since $h( r_\mathrm{S} z)$ is a univalent mapping from 
the unit disk to the Siegel disk, we can 
use the area formula \cite{Rudin87}
\begin{equation}\label{areaformula} 
{\rm Area} = \pi \sum_{k = 1}^\infty k |h_k|^2 
= 
\pi \sum_{k = 1}^\infty k r_\mathrm{S}^{-2k} |\hat \chi_{-k}|^2 \ .
\end{equation} 

For polynomials, the Siegel disk is bounded
so  that the sum in \eqref{areaformula} 
is finite. This is compatible with the 
observation \eqref{eq:decay-modulus}, 
but it shows that the bound cannot be saturated very often.

In Table~\ref{table:areas} 
\begin{table*}
\caption{Areas of the Siegel disks for maps of the form 
        $f_{d,\langle k^\infty\rangle,1+3i}$ (see~\eqref{eq:f-def})
        with different rotation numbers $\sigma$ and orders $d$ 
        of the critical point.\label{table:areas}
}
\begin{ruledtabular}
\begin{tabular}{rlllll}
\mbox{$d$} & \mbox{$\langle 1^\infty\rangle$}  & \mbox{$\langle 2^\infty\rangle$} & \mbox{$\langle 3^\infty\rangle$} & \mbox{$\langle 4^\infty\rangle$} & \mbox{$\langle 5^\infty\rangle$} \\
\hline
 1  & 1.3603361 & 1.3586530 & 1.3611085 & 1.3652030 & 1.3693337 \\
 2  & 0.895659  & 0.893442  & 0.893605  & 0.89408   & 0.89367   \\
 3  & 0.65986   & 0.65766   & 0.65664   & 0.6553    & 0.65308   \\
 4  & 0.5190    & 0.5170    & 0.51550   & 0.5133    & 0.5104    \\
 5  & 0.4262    & 0.4244    & 0.4226    & 0.4200    & 0.4170    \\
 6  & 0.3607    & 0.3591    & 0.3573    & 0.354     & 0.3515    \\
10  & 0.2214    & 0.2203    & 0.2187    & 0.2163    & 0.2137    \\
15  & 0.148     & 0.147     & 0.146     & 0.144     & 0.1421    \\
20  & 0.111     & 0.110     & 0.109     & 0.107     & 0.106     \\
\end{tabular}
\end{ruledtabular}
\end{table*}
we give the values of the areas of the Siegel disks 
of the map $f_{d,\sigma,\beta}$ \eqref{eq:f-def} 
(we believe that the error does not exceed 2 
in the last digit).  
Since the series \eqref{areaformula} converges slowly, 
we computed the partial sums of the first $Q_n$ terms 
in \eqref{areaformula}, where $Q_n$ are the denominators 
of the best rational approximants, 
$P_m/Q_m = \langle k^m\rangle$, 
to the rotation number $\sigma=\langle k^\infty\rangle$ 
(cf.~\eqref{eq:fibonacci}), 
and then performed Aitken extrapolation on these values.  
Because of the repeating ``periodic windows'' 
in the Fourier spectra 
(shown in Figures \ref{fig:spectrum_1_inf}, \ref{fig:spectrum_5_inf}, 
\ref{fig:spectra_1_inf_and_4_inf_straight}, 
\ref{fig:spectrum_2_inf_crit_1_5_20_reduced}), 
these partial sums tend to the area of the Siegel disk geometrically, 
and Aitken extrapolation gives good results.  
In our computations we used $2^{23}\approx 8\times 10^6$ 
Fourier coefficients of~$\chi$.  

Clearly, the area of a Siegel disk depends on the particular 
choice of the map~$f$, i.e., is non-universal.  
Perhaps the only universal characteristic that can be extracted 
is the rate of convergence in the Aitken extrapolation, 
but we have not studied this problem in detail.

\section{An upper bound of the regularity of 
the conjugacy\label{sec:bound-regularity}} 

As pointed out in 
\cite[Section~8.2]{LlaveP02}, one 
can find upper bounds for the regularity in
terms of the  scaling exponents.

Recall that $\chi:\bbT\to\partial\siegeldisk$ 
conjugates (the restriction to $\partial\siegeldisk$ of) 
the map $f:\partial\siegeldisk\to\partial\siegeldisk$ 
to the rigid rotation 
$r_\sigma:\bbT\to\bbT:t\mapsto t+\sigma$ 
(where $\sigma$ is the rotation number of~$f$), 
namely, 
$\chi\circ r_\sigma = f\circ \chi$.  
Let us consider only rotation numbers of the form 
$\sigma=\langle k^\infty \rangle$, 
and let the natural numbers $Q_m$ and 
the scaling exponent $\alpha$ be defined 
by \eqref{eq:fibonacci} and~\eqref{eq:scaling-exp}.  
Then closest returns of the iterates of $0\in\bbT$ of $r_\sigma$ 
and the iterates of $c\in\partial\siegeldisk$ 
to the starting points $0$ and $c$, respectively, 
are governed by the scaling relations 
\begin{eqnarray*}
r_\sigma^{Q_m}(0) &=& C_1 \sigma^m + o(\sigma^m) \ , \\
|f^{Q_m}(c)-c| &=& C_2 |\alpha|^m + o(|\alpha|^m) \ .
\end{eqnarray*}
So that we obtain 
that $h( C_1 \sigma^m) \approx C_2 |\alpha|^m$. 

This is impossible if $h$ is $C^\kappa$ with 
\begin{equation}
\label{eq:max-regul}
\regul > \regul_\mathrm{max} \equiv 
\log |\alpha|/ \log \sigma
\end{equation}
and the right-hand side of  \eqref{eq:max-regul} is not an integer. 

In  Table~\ref{table:regularity-scaling}, 
%%%%
\begin{table}
\caption{Upper bounds on  the H\"older regularity 
        $\regul_\mathrm{max}$ \eqref{eq:max-regul} 
        for maps of the form $f_{d,\sigma,\beta}$ \eqref{eq:f-def} 
        with rotation numbers 
        $\sigma=\langle k^\infty\rangle$ ($k=1,2,3,4,5$) 
        and critical point of order~$d$.  
        \label{table:regularity-scaling}
}
\begin{ruledtabular}
\begin{tabular}{rccccc}
$d$ & $\langle 1^\infty\rangle$ & $\langle 2^\infty\rangle$ & $\langle 3^\infty\rangle$ & $\langle 4^\infty\rangle$ & $\langle 5^\infty\rangle$\\
\hline
 1 &  0.6203 & 0.6159 & 0.6064 & 0.5933 & 0.5783 \\
 2 &  0.4324 & 0.4276 & 0.4175 & 0.4039 & 0.3890 \\
 3 &  0.3293 & 0.3252 & 0.3164 & 0.3047 & 0.2922 \\
 4 &  0.2653 & 0.2618 & 0.2543 & 0.2443 & 0.2338 \\
 5 &  0.2219 & 0.2189 & 0.2124 & 0.2039 & 0.1949 \\
 6 &  0.1906 & 0.1880 & 0.1823 & 0.1749 & 0.1670 \\
10 &  0.1218 & 0.1200 & 0.1163 & 0.1114 & 0.1063 \\
15 &  0.0838 & 0.0826 & 0.0800 & 0.0766 & 0.0731 \\
20 &  0.0639 & 0.0630 & 0.0610 & 0.0584 & 0.0557 \\
\end{tabular}
\end{ruledtabular}
\end{table}
%%%%
we give the values of the upper bound to the  regularity 
$\regul_\mathrm{max}$ \eqref{eq:max-regul} 
for maps with rotation numbers and order of the critical point.  
To compute these values, 
we used the values for $|\alpha|$ 
from Table~\ref{table:exponents} 
(and the exact values for the rotation numbers).  
We have kept only four digits of accuracy, 
although the error of these numbers is smaller 
(their relative error is the same as the relative error 
of the values of~$|\alpha|$).  

Clearly, within the numerical error, 
the values of the regularities from Table~\ref{table:regularity} 
(obtained from applying the CLP method) 
are equal to the upper bounds on the regularity 
from Table~\ref{table:regularity-scaling} 
(obtained from the scaling exponents).

This is in contrast with the results in 
\cite{LlaveP02}, where similar bounds based on scaling 
were found to be saturated by some conjugacies
but not by the inverse conjugacy.

For maps with highly critical points, 
if \eqref{eq:max-regul} holds, 
then the asymptotic behavior of $|\alpha|$ \eqref{eq:A} 
implies that the asymptotic behavior 
of the limit on the regularity becomes 
\begin{equation}
\regul_{\mathrm{max},\sigma} 
\sim 
\frac{\log \left|1-\frac{A_\sigma}{d}\right|}{\log \sigma} 
\approx  
\frac{A_\sigma}{|\log|\sigma||} \, \frac{1}{d} 
\qquad \mbox{ as } \ d\to\infty \ .
\end{equation}

%%%%%%%%%%%%%%%%%%%%%%%%%%%%%%%%%%%%%%%%%%%%%%%%%%%%%%%%%%%%

\section{Conclusions\label{sec:discussion}} 

We have considered Siegel disks
of polynomials with some quadratic fields and 
with different degrees of critical points. 

We have made extended precision calculations of 
scaling exponents and a parameterization of the boundary. 
This allows us to compute the regularity of the boundary 
using methods of harmonic analysis. 

The regularity of the boundary seems to be universal, 
depend only on the tail of the continued fraction expansion 
and saturate some easy bounds in terms of the continued fraction 
expansion. 

We have identified several regularities of the Fourier series of 
the conjugacy. Namely, it seems that, irrespective of the rotation 
number, we have $0 < \limsup |k\hat \chi_k| <\infty$. There seems to be 
a regular distribution of the phases of the Fourier coefficients 
which follows a Gaussian law. 

We have also extended the results of \cite{Osbaldestin92} on 
the dependence of scaling exponents on the degree of the 
critical point to higher degrees and to other rotation numbers.

%\section*{Acknowledgments}

\begin{acknowledgments}
The work of both authors has been partially supported by NSF grants. 
We thank Amit Apte, Lukas Geyer, Arturo Olvera for useful discussions.
This work would not have been possible without access to 
excellent publicly available software for 
compiler, graphics, extended precision, especially GMP~\cite{GMP}.  
We thank the Mathematics Department at University of Texas at Austin for 
the use of computer facilities.
\end{acknowledgments}

%\begin{acknowledgments}
%\end{acknowledgments}

%\bibliography{llave99,new,new1}
\def\cprime{$'$} \def\cprime{$'$} \def\cprime{$'$}

\end{document}